\documentclass{amsart}
\usepackage{amsmath,amsthm,amssymb,mathrsfs,mathtools,MnSymbol,IMjournal}
\usepackage{algorithmic}
\usepackage[ruled]{algorithm}

\usepackage{color}
\definecolor{dkr}{rgb}{.6,0,0}
\definecolor{dkg}{rgb}{0,.6,0}
\definecolor{dkb}{rgb}{0,0,.8}
\definecolor{dko}{rgb}{.6,0,.2}

\usepackage[colorlinks=true, 
  citecolor=dkg, linkcolor=dkr, urlcolor=dkb, %
  bookmarks=false, pdfpagelabels=false, 
  hypertexnames=false, naturalnames=true, 
  pagebackref=false]{hyperref}

\renewcommand{\vec}{{\mathrm{vec}}}
\newcommand{\diag}{{\mathrm{diag}}}
\newcommand{\trp}{{\mathsf{T}}}
\newcommand{\fro}{{\mathsf{F}}}
\newcommand{\kro}{{\mathsf{K}}}

\renewcommand{\hat}{\widehat}

\newcommand{\FDL}{\Delta_{\mathsf{FD}}}
\newcommand{\LapDD}{L_{\mathsf{2D}}}
\newcommand{\LapDDD}{L_{\mathsf{3D}}}

\newcommand{\Lap}[2]{L_{{#1},{#2}}}
\newcommand{\LapP}[1]{\Lap{#1}{\mathsf{P}}}
\newcommand{\LapD}[1]{\Lap{#1}{\mathsf{D}}}
\newcommand{\LapN}[1]{\Lap{#1}{\mathsf{N}}}
\newcommand{\LapDN}[1]{\Lap{#1}{\mathsf{DN}}}
\newcommand{\LapND}[1]{\Lap{#1}{\mathsf{ND}}}

\begin{document}
\newtheorem{theorem}{Theorem}[section]
\newtheorem{lemma}[theorem]{Lemma}
\newtheorem{remark}[theorem]{Remark}
\newtheorem{problem}[theorem]{Problem}
\newtheorem{experiment}[theorem]{Experiment}
\numberwithin{equation}{section}

\title{Tensor-structured PCG for finite difference solver of domain patterns in ferroelectric material}

\author{\|V\v{e}nceslav |Chumchal|, Liberec,
        \|Pavel |Marton|, Liberec \& Praha,
        \|Martin |Ple\v{s}inger|, Liberec,
        \|Martina |\v{S}im\r{u}nkov\'{a}|, Liberec}



\abstract 
This paper presents a case study of application of the preconditioned method of conjugate gradients (CG)
on a problem with operator resembling the structure of sum of Kronecker products. 
In particular, we are solving the Poisson's equation on a sample of homogeneous 
isotropic ferroelectric material of cuboid shape, where the Laplacian is discretized by 
finite difference. We present several preconditioners that fits the Kronecker structure 
and thus can be efficiently implemented and applied. Preconditioner based on the Moore--Penrose 
pseudoinverse is extremely efficient for this particular problem, and also applicable 
(if we are able to store the dense right-hand side of our problem).
We briefly analyze the computational cost of the method and individual preconditioners,
and illustrate effectiveness of the chosen one by numerical experiments. 

Although we describe our method as {\em preconditioned CG with pseudoinverse-based preconditioner},
it can also be seen as {\em pseudoinverse-based direct solver with iterative refinement by CG iteration}.

This work is motivated by real application, the method was already implemented in C/C++ code
{\sc Ferrodo2} and first results were published in Physical Review B 107(9) (2023), paper id 094102.
\endabstract




\section{Introduction}\label{sec:intro}

We are interested in the computation of the electric field induced by the given 
distribution of electric charge density $\rho$ within a dielectric material. 
For this, we use Poisson's equation
\[
  -\Delta \upsilon = \eta
\]
which links electric potential and charge density ($\Delta$ is the Laplace operator). 
In particular, $\eta$ represents the rescaled density of electrostatic charge ($\eta =\rho/\epsilon$, 
where $\epsilon$ is the electric permitivity of the material, which is considered here homogeneous and scalar), 
$\upsilon$ denotes the electric potential induced by $\rho$.

The Poisson's equation directly follows from the third Maxwell's equation ($\nabla\cdot\delta = \eta$), 
relation between electric induction $\delta$ and the target electric field $\varepsilon$ ($\delta=\epsilon \varepsilon$), 
and from the relation between electric potential $\upsilon$ and electric field 
\[
  \varepsilon = -\nabla\cdot\upsilon.
\]
The existence of the potential is provided by the second Maxwell's equation ($\nabla\times\varepsilon=0$) for the electrostatics. 
Once the potential $\upsilon$ is recovered from Poisson's equation, the sought vector field --- the electric field $\varepsilon$ can 
then be easily obtained from $\upsilon$ using the second equation. For more details about the electrostatic considerations, 
see \cite{maxwell} or, e.g., \cite{feynman2}. Also, it is fair to note here, that we do not follow the standard notation of 
the physical quantities to avoid possible collisions with also very standard mathematical notation later on.

The ultimate goal of a wider work than presented here is to contribute to understanding the so-called ferroelectric domain 
structures that emerge within the ferroelectric materials, which are utilized in a broad range of applications in actuation, 
sensing, and beyond. The domain structures are known to have a critical impact on material's properties, and at the same time, 
they are very difficult to determine experimentally, particularly for nano-size domains. Therefore, numerical simulations of domain 
structures is a necessary prerequisite for their understanding and control. An important aspect of these simulations is determining 
how these domains are influenced by charged defects/dopants within the ferroelectric materials. In this paper, however, we focus 
only on the numerical estimation of the electric potential $\upsilon$, i.e., on approximate solution of the discretized Poisson's equation.

\subsection{Problem description and its discretization}

Although our goal is to solve real-world-like, so three-dimensional (3D) problems,
we start with simpler model problems that are only two-dimensional (2D). We are
interested only in rectangular 2D or cuboid 3D domains. This choice is suitable 
from the computational point of view because it imposes Kronecker structure to 
the problem. However, it is also relevant from the practical, application point 
of view --- we often deal with periodic boundary condition of the sample, where
the cuboid can be seen as a single cell in larger periodic structure, which, e.g., 
opens door to Fourier approach. 

On each pair of the opposite sides, or faces of the rectangle, or cuboid, respectively,
there can be one of the following four types of boundary conditions (BC):
\begin{itemize}\itemsep 0pt
\item Periodic BC (e.g., the electric potential on the left side of rectangle
is the same as on the right side; it reflects periodicity of microscopic structure of the material);
\item Dirichlet BC (which prescribes electric potential on the boundary);
\item Neumann BC (which prescribes derivative of electric potential, i.e., the electric field over the boundary); or
\item Mixed BC, in particular Dirichlet BC on one side (or face) and Neumann BC on the opposite.
\end{itemize}
Notice that for some combinations of boundary conditions ---
in particular periodic or Neumann in both (all three) directions ---
the Laplacian is not invertible, i.e., the solution is not unique,
it is given up to a constant electric potential (notice, this does not affect
its gradient and thus the electric field).

Due to the simple geometry, we consider the finite difference (FD) discretization.
Since such discretization of the Laplace operator in two (or three) mutually
orthogonal directions is trivial, the only two pieces of information we have
in hands are: The discretized right-hand side $\eta$, i.e., a matrix $H\in\mathbb{R}^{n\times q}$
or its three-way analogue --- a tensor $\mathcal{H}\in\mathbb{R}^{n\times q\times t}$,
and the structure of BCs. We are looking for discretized electric potential
$\upsilon$, i.e., an unknown matrix $U\in\mathbb{R}^{n\times q}$ or an unknown
tensor $\mathcal{U}\in\mathbb{R}^{n\times q\times t}$.

\subsection{Notation and structure of the work}

Throughout the paper we try to keep the convention that vectors are denoted by lowercase
italic letters, whereas matrices by uppercase italic letters, and tensors by uppercase
calligraphic letters; entries of these objects are also denoted by lowercase italic letter,
for example $(A)_{i,j}=a_{i,j}$, or $(\mathcal{A})_{i,j,k}=a_{i,j,k}$. In particular $I$
and $I_n$ denotes the identity matrix of suitable or given size. 
Further, $1$, $1_n$, $1_{n,q}$, and $1_{n,q,t}$ denotes number one, vector of
ones, matrix of ones, or tensor of ones of suitable or given size; similarly $0$, $0_n$,
$0_{n,q}$, and $0_{n,q,t}$. 
By $\otimes$ we denote Kronecker product of matrices,
by $\odot$ Hadamard (entry-wise) product of matrices,
by $\times_\ell$ the $\ell$-mode matrix-tensor product.
As usual, $A^{-1}$ denotes inverse of $A$ (if it exists),
$A^\dagger$ the Moore--Penrose pseudoinverse,
and $A^{\odot-1}$, $A^{\odot\dagger}$ their Hadamard (entry-wise) versions.
By $\lambda$ or $\lambda(A)$ we denote the eigenvalue
of $A$, by $\lambda_k$ or $\lambda_k(A)$ we denote the $k$th largest one (including
multiplicities), and $\sigma(A)$ denotes the spectrum (set of all eigenvalues) of $A$.
Finally $\|v\|_2$ stands for the Euclidean (i.e., two) norm of vectoR $v$ and $\|A\|_\fro$
for the Frobenius norm of matrix $A$.

The text is organized as follows: After the introduction in Section~\ref{sec:intro},
we recapitulate some basic concepts and notation related to the tensor manipulation
in Section~\ref{sec:tensorintro}. Section~\ref{sec:algsystem} describes the basic
properties and structure of linear algebraic system originated in discrete Laplacian.
Section~\ref{sec:PCG} discusses application of the preconditioned method of conjugate
gradients. The obtained results are presented in Section~\ref{sec:numexp}.
Section~\ref{sec:conclusions} concludes the paper.

\section{Basic tensor concepts and notation}\label{sec:tensorintro}

We are dealing with 2D (in the simplified model), or later also 3D data arranged
into regular grid of rectangular, or cuboid shape (thanks to the cuboid shape
of the analyzed sample of material and regular finite difference discretization).
These data can be naturally written in the form of a matrix, or its three-way
analogy, respectively. We call the latter simply a {\em tensor} as it is usual
in the matrix computations community. It would be useful here to introduce
the tensor-related notation used throughout the text. All the concepts below can
be found, e.g., in the review paper of Kolda and Bader \cite{KolBad09}
and in the references therein.

\subsection{Vectorization}\label{ssec:vectorization}

Let $U\in\mathbb{R}^{n\times q}$ be a matrix with
columns $u_j$ and entries $u_{i,j}$, where $i=1,\ldots,n$, $j=1\ldots,q$. By {\em vectorization} we understand
mapping of $U$ to one long column-vector
\[
  \begin{split}
  \vec(U)&=[u_1^\trp,u_2^\trp\ldots,u_q^\trp]^\trp \\
         &=\big[
          [u_{1,1},u_{2,1},\ldots,u_{n,1}]\,,\,
          [u_{1,2},u_{2,2},\ldots,u_{n,2}]\,,\,
          \ldots\,,\,
          [u_{1,q},u_{2,q},\ldots,u_{n,q}]
        \big]^\trp\in\mathbb{R}^{nq}
  \end{split}
\]
that collects all the columns of $U$, one beneath other. Equivalently, $\vec(U)$ contains all
the entries $u_{i,j}$ of $U$ sorted in the inverse lexicographical order w.r.t. their multiindices $(i,j)$; see for example \cite{Turkington13}.

The equivalent formulation can be directly generalized for a tensor $\mathcal{U}\in\mathbb{R}^{n\times q\times t}$
with entries $u_{i,j,k}$ described by multiindices $(i,j,k)$, i.e.,
\[
  \begin{split}
  \vec(\mathcal{U})=\big[
    &[u_{1,1,1},\ldots,u_{n,1,1}]\,,\,
    [u_{1,2,1},\ldots,u_{n,2,1}]\,,\,
    \ldots\,,\,
    [u_{1,q,1},\ldots,u_{n,q,1}]\,,\,
    \\
    &[u_{1,1,2},\ldots,u_{n,1,2}]\,,\,
    [u_{1,2,2},\ldots,u_{n,2,2}]\,,\,
    \ldots\,,\,
    [u_{1,q,2},\ldots,u_{n,q,2}]\,,\,
    \\
    &\ldots\,,\,
    \\
    &[u_{1,1,t},\ldots,u_{n,1,t}]\,,\,
    [u_{1,2,t},\ldots,u_{n,2,t}]\,,\,
    \ldots\,,\,
    [u_{1,q,t},\ldots,u_{n,q,t}]
   \big]^\trp\in\mathbb{R}^{nqt}.
  \end{split}
\]
(The inner brackets are there only for better legibility.)

\subsection{Matrix-tensor product (MT)}\label{ssec:MTproduct}

First recall the {\em Kronecker product} of two matrices
$A\in\mathbb{R}^{m\times n}$ and $B\in\mathbb{R}^{p\times q}$, i.e.,
\[
  B\otimes A = \left[\begin{array}{cccc}
    b_{1,1}A & b_{1,2}A & \cdots & b_{1,q}A \\
    b_{2,1}A & b_{2,2}A & \cdots & b_{2,q}A \\
    \vdots   & \vdots   & \ddots & \vdots   \\
    b_{p,1}A & b_{p,2}A & \cdots & b_{p,q}A
  \end{array}\right]\in\mathbb{R}^{(mp)\times(nq)}.
\]
It has a lot of interesting properties (see \cite{b:fiedler} or \cite{Turkington13}),
and a lot of relations to the standard matrix summation and multiplication.
If $A_1$ and $A_2$ are matrices of the same dimensions, and also $B_1$ and $B_2$, then
\begin{equation}\label{eq:kronPsum}
  (B_1+B_2)\otimes (A_1+A_2) = (B_1\otimes A_1)+(B_1\otimes A_2)+(B_2\otimes A_1)+(B_2\otimes A_2).
\end{equation}
On the other hand, if products $A_1A_2$ and $B_1B_2$ exist, then
\begin{equation}\label{eq:kronPprod}
  (B_1B_2)\otimes (A_1A_2) = (B_1\otimes A_1)(B_2\otimes A_2).
\end{equation}
It relates also to vectorization; if $A\in\mathbb{R}^{m\times n}$,  $U\in\mathbb{R}^{n\times p}$,
and $B\in\mathbb{R}^{p\times q}$, then
\begin{equation}\label{eq:kronPvec}
  \vec(AUB^\trp) = (B\otimes A)\,\vec(U).
\end{equation}
Notice that in $AUB^\trp$ the matrix $A$ multiplies columns of $U$ whereas
$B$ multiplies rows of $U$ (both treated as column-vectors). In terms of entries:
columns of $U$ are collections of entries $u_{i,j}$ with variable first index $i$
and all the other indices (here only the second $j$) is fixed; similarly rows are
collections of entries $u_{i,j}$ with variable second index $j$, whereas all the
other indices (here only the first $i$) is fixed.

The entry-wise point of view allows us again to generalize several concepts to tensors.
Any collection of all entries $u_{i,j,k}$ with one particular, let say $\ell$th ($\ell\in\{1,2,3\}$)
index variable (e.g., for $\ell=1$, $i=1,2,\ldots,n$) and all the other indices fixed (e.g., $(j,k)=(3,7)$)
is called the {\em $\ell$-mode fiber} (in our case $(3,7)$th $1$-mode fiber).
The {\em $\ell$-mode matrix-tensor product} $M\times_\ell\mathcal{U}$ is then
multiplication of $M$ with all the {\em $\ell$-mode fibres} (treated
as column-vectors). In our case where
$A\in\mathbb{R}^{m\times n}$,
$B\in\mathbb{R}^{p\times q}$,
$C\in\mathbb{R}^{s\times t}$, and
$\mathcal{U}\in\mathbb{R}^{n\times q\times t}$, we get
\[
  \begin{split}
  A\times_1\mathcal{U}\in\mathbb{R}^{m\times q\times t}, \qquad
  &\vec(A\times_1\mathcal{U}) = (I_t\otimes I_q\otimes A  )\,\vec(\mathcal{U}), \\
  B\times_2\mathcal{U}\in\mathbb{R}^{n\times p\times t}, \qquad
  &\vec(B\times_2\mathcal{U}) = (I_t\otimes B  \otimes I_n)\,\vec(\mathcal{U}), \\
  C\times_3\mathcal{U}\in\mathbb{R}^{n\times q\times s},  \qquad
  &\vec(C\times_3\mathcal{U}) = (C  \otimes I_q\otimes I_n)\,\vec(\mathcal{U}). \\
  \end{split}
\]
Multiplication of a tensor by two matrices in two different modes is independent
on the order of individual operations (which follows directly from the associativity
of matrix multiplication). Therefore, we simplify the notation a bit in the case of
{\em general linear transformation} of the tensor
\[
  (A,B,C\,|\,\mathcal{U}) = A\times_1(B\times_2(C\times_3\mathcal{U}))\in\mathbb{R}^{m\times p\times s}.
\]
Similarly as before,
\[
  \vec((A,B,C\,|\,\mathcal{U})) = (C\otimes B\otimes A)\,\vec(\mathcal{U}).
\]
For deeper insight into tensor manipulation we recommend to study \cite{KolBad09}.

Finally, notice that in our case the transforming matrices $A$, $B$, and $C$ will
always be square, symmetric, and positive semi-definite.

\section{Description of the algebraic system}\label{sec:algsystem}

As already mentioned, we are focused on solving the discretized Poisson's equation
\[
  -\Delta \upsilon = \eta;
\]
see \cite{evans}.
We use regular 2D or 3D finite difference discretization with the mesh oriented
in accordance with the sample of material. We moreover assume, for simplicity,
that the material is homogenous and isotropous, i.e., the material coefficients
does not change along the discretization and in different orientations. Therefore
we have very simple discretization of the (minus) Laplace operator.

\subsection{Discretization of 1D Laplace operator}

Recall that in one-dimensional (1D) case, the finite difference discretization of 
the {\em minus} Laplacian takes exceptionally simple form
\begin{equation}\label{eq:lap1d}
  L_{n}(\alpha,\beta,\gamma) = \left[\begin{array}{cccccc}
  \alpha & -1 & & & & \gamma \\
  -1 & 2 & -1 \\
  & -1 & 2 & \ddots \\
  & & \ddots & \ddots & -1 \\
  & & & -1 & 2 & -1 \\
  \gamma & & & & -1 & \beta
  \end{array}\right]
  \in\mathbb{R}^{n\times n}, \qquad n \geq 3,
\end{equation}
see, e.g., \cite{b:fiedler}, \cite{b:varga}.
The matrix is obviously symmetric $L_{n}(\alpha,\beta,\gamma) = (L_{n}(\alpha,\beta,\gamma))^\trp$.
Values of entries in its corners depend on BCs, in particular:
\begin{itemize}\itemsep 0pt
\item Periodic BC enforces $(\alpha,\beta,\gamma) = (2,2,-1)$;
\item Dirichlet BC enforces $(\alpha,\beta,\gamma) = (2,2,0)$;
\item Neumann BC enforces $(\alpha,\beta,\gamma) = (1,1,0)$; and
\item Mixed BC enforces $(\alpha,\beta,\gamma) = (2,1,0)$ or $(1,2,0)$.
\end{itemize}
If it is not so important, then we will denote the matrix simply $L_n$ or even $L$.
On the other hand, if we need to distinguish the particular BC/matrix structure, we use
\[
  \begin{array}{c}
  \LapP{n}=L_n(2,2,-1), \quad
  \LapD{n}=L_n(2,2,0), \quad
  \LapN{n}=L_n(1,1,0),
  \\[2mm]
  \LapDN{n}=L_n(2,1,0), \quad \text{and} \quad
  \LapND{n}=L_n(1,2,0).
  \end{array}
\]
Notice, that the BC affects in general also the right-hand side of the algebraic
systems; see Section \ref{ssec:modifyRHS}.

\subsection{Spectral properties of $L_n$ matrices}

Here we briefly mention the basic spectral properties (eigenvalues and
eigenpairs) of matrices $L_n$. We listed them without proofs,
but these results are commonly known and can be easily verified (simply by
multiplication the matrices by respective vectors and employing basic trigonometry);
see for example \cite{kolarova}.


\begin{lemma}[Spectral properties of $\LapP{n}$]\label{lem:spectrumP}
Let $n\in\mathbb{N}$, $n\geq 3$. Consider $\ell\in\mathbb{N}$ such that either $n = 2\ell-1$ or $n = 2\ell$.
Then
\[
  \lambda_{1}=0, \qquad
  \lambda_{2k} = \lambda_{2k+1} = 2-2\cos\bigg(\frac{2k\pi}{n}\bigg), \quad k=1,\ldots,\ell-1, \qquad
  \lambda_{2\ell}=4,
\]
represent the eigenvalues and
\[
  \begin{split}
  v_{1}     &= [1,1,\ldots,1]^\trp, \\
  v_{2k}    &= \left[\cos\bigg(\frac{2jk\pi}{n}\bigg)\right]_{j=1}^n, \quad
  v_{2k+1}   = \left[\sin\bigg(\frac{2jk\pi}{n}\bigg)\right]_{j=1}^n, \quad k=1,\ldots,\ell-1, \\
  v_{2\ell} &= [1,-1,\ldots,1,-1]^\trp.
  \end{split}
\]
the eigenvectors of $\LapP{n}$ (except of the last pair in the odd $n$ case).
\end{lemma}


\begin{lemma}[Spectral properties of $\LapD{n}$]\label{lem:spectrumD}
Let $n\in\mathbb{N}$, $n\geq 3$. Then
\[
  \lambda_{k} = 2-2\cos\bigg(\frac{k\pi}{n+1}\bigg), \quad
  v_{k} = \left[\sin\bigg(\frac{jk\pi}{n+1}\bigg)\right]_{j=1}^n,
  \quad k=1,\ldots,n,
\]
represent the eigenvalues and the eigenvectors of $\LapD{n}$.
\end{lemma}


\begin{lemma}[Spectral properties of $\LapN{n}$]\label{lem:spectrumN}
Let $n\in\mathbb{N}$, $n\geq 3$. Then
\[
  \lambda_{k} = 2-2\cos\bigg(\frac{(k-1)\pi}{n}\bigg), \quad
  v_{k} = \left[\cos\bigg(\frac{(j-\frac12)(k-1)\pi}{n}\bigg)\right]_{j=1}^n,
  \quad k=1,\ldots,n,
\]
represent the eigenvalues and the eigenvectors of $\LapN{n}$.
\end{lemma}


\begin{lemma}[Spectral properties of $\LapDN{n}$ and $\LapND{n}$]\label{lem:spectrumM}
Let $n\in\mathbb{N}$, $n\geq 3$. Then
\[
  \lambda_{k} = 2-2\cos\bigg(\frac{(2k-1)\pi}{2n+1}\bigg), \quad
  v_{k} = \left\{\begin{array}{l}
  \displaystyle \left[\sin\bigg(\frac{j(2k-1)\pi}{2n+1} - k\pi\bigg)\right]_{j=1}^n \\[2mm]
  \displaystyle \left[\cos\bigg(\frac{j(2k-1)\pi}{2n+1} - \frac{(2k-1)\pi}{2(2n+1)}\bigg)\right]_{j=1}^n
  \end{array}\right\},
\]
$k=1,\ldots,n$, represent the eigenvalues and the eigenvectors of $\LapDN{n}$ and $\LapND{n}$, respectively.
\end{lemma}


The goal of this list is to point out several things. First, we see, that all the
spectral information of all the matrices is available in analytical form. Second,
we see, that essentially all the eigenvalues are of the form $2-2\cos(f(j)\pi/g(n))$,
where $f(x)$ and $g(x)$ are linear functions, $f(x)$ is nonnegative and $g(x)$ positive
on the respective domain and $f(j)<g(n)$. Thus, all the eigenvaues of all the matrices
are nonnegative from interval $[0,4]$. Third, actually most of the eigenvalues are
positive --- there are only two exceptions
\[
  \lambda_1(\LapP{n}) = \lambda_1(\LapN{n}) = 0.
\]
Consequently, matrices $\LapD{n}$, $\LapDN{n}$, $\LapND{n}$ are always symmetric
positive definite, whereas $\LapP{n}$ and $\LapN{n}$ are always only symmetric
positive semi-defnite, i.e., not invertible, $\det(\LapP{n})=\det(\LapN{n})=0$.
Finally, notice that matrices $\LapDN{n}$ and $\LapND{n}$ in Lemma \ref{lem:spectrumM}
are essentially the same, just upside-down \& left-to-right rearranged. Therefore,
also their eigenvectors are essentially the same (only upside-down rearranged,
i.e., reindexed by $j\,\longmapsto\,n+1-j$ and after few goniometric manipulations
multiplied by $-1$).

\subsection{Discretization of 2D \& 3D Laplacians and their spectra}

In the 1D case, our equation $-\Delta\upsilon=\eta$ is by {\em finite differences} (FD)
projected to the discrete space as
\[
  -\FDL(u) \equiv L_n u = h.
\]
In 2D and 3D cases is situation very similar. Thanks to the
discretization independent in individual directions and thanks to simple scalar
homogenous properties of the material, we get
\begin{equation}\label{eq:ourproblem}
  -\FDL(U) \equiv L_n U + U L_q^\trp = H,
  \quad
  -\FDL(\mathcal{U}) \equiv L_n\times_1\mathcal{U} + L_q\times_2\mathcal{U} + L_t\times_3\mathcal{U} = \mathcal{H},
\end{equation}
respectively; see, e.g., \cite{b:fiedler}, \cite{b:varga}. By vectorization of these equations
we get them in the standard form with single matrix, vector of unknowns and
the right-hand side vector
\begin{equation}\label{eq:ourproblem2}
  \begin{split}
  \Big(\overbrace{(I_q\otimes L_n)+(L_q\otimes I_n)}^{\displaystyle \LapDD\in\mathbb{R}^{(nq)\times(nq)}}\Big)\,\vec(U) &= \vec(H), \\
  \Big(\underbrace{(I_t\otimes I_q\otimes L_n)+(I_t\otimes L_q\otimes I_n)+(L_t\otimes I_q\otimes I_n)}_{\displaystyle \LapDDD\in\mathbb{R}^{(nqt)\times(nqt)}}\Big)\,\vec(\mathcal{U}) &= \vec(\mathcal{H}).
  \end{split}
\end{equation}
It is clear that these matrices are unnecessarily big, since all the useful
information is carried only by the two, or three bits of information ---
the individual $L_\ell$ matrices. These big matrices are cannot be used for practical
computation, but they can provide us a lot of useful information.

Since our matrix $L_\ell$ is symmetric, we can always write it in the form of spectral
decomposition, i.e., as a product
\begin{equation}\label{eq:eigL2D1}
  L_\ell = V_\ell\Lambda_\ell V_\ell^{-1}.
\end{equation}
Here $\Lambda_\ell$ is a diagonal matrix with diagonal entries equal to eigenvalues of $L_\ell$,
and $V_\ell$ is an invertible matrix with columns equal to eigenvector of $L_\ell$ (both given by
Lemmas~\ref{lem:spectrumP}--\ref{lem:spectrumM}; notice that the eigenvectors
are not normalized there). Employing (\ref{eq:kronPprod}), we get
\begin{equation}\label{eq:eigL2D2}
  \begin{split}
  \LapDD
  &= \Big((V_qV_q^{-1}\otimes V_n\Lambda_nV_n^{-1})+(V_q\Lambda_qV_q^{-1}\otimes V_nV_n^{-1})\Big) \\
  &= (V_q\otimes V_n)\,\Big((I_q\otimes \Lambda_n)+(\Lambda_q\otimes I_n)\Big)\,(V_q\otimes V_n)^{-1},
  \end{split}
\end{equation}
which is the spectral decomposition, and thus
\[
  \sigma(\LapDD) = \big\{
  \lambda_i(L_n)+\lambda_j(L_q)\;:\quad i=1,\ldots,n, \; j=1,\ldots q
  \big\}
\]
is the spectrum of the discretized 2D minus Laplacian. Similarly,
\[
  \sigma(\LapDDD) = \big\{
  \lambda_i(L_n)+\lambda_j(L_q)+\lambda_k(L_t)\;:\quad i=1,\ldots,n, \; j=1,\ldots q, \; k=1,\ldots,t
  \big\}
\]
is the spectrum of the discretized 3D minus Laplacian.

Finally, we see that the big matrices $\LapDD$ and $\LapDDD$ in (\ref{eq:ourproblem2})
are always symmetric positive semi-definite. Moreover, they are invertible (and thus
positive definite and thus (\ref{eq:ourproblem}) uniquely solvable) if and only if
at least one of the matrices $L_\ell$ is invertible.

\subsection{Modifications of the right-hand side enforced by BCs}\label{ssec:modifyRHS}

Since the solving process of $\Delta \upsilon = \eta$ starts with the discretized
right-hand side (obtained from the outside) and the sets of BCs, we need to ensure the
compatibility of these two pieces of information. In particular the BCs will govern
the structure of right-hand side. It will be done in a very standard way, thus we
do only a few comments.

\subsubsection{Orthogonality to null-space}\label{sssec:RHS-nullspace}

First, recall that discretized 1D Laplacians $\LapP{n}$ and $\LapN{n}$ are singular,
both with simple zero eigenvalue. Similarly, discretized 2D and 3D Laplacians
(\ref{eq:ourproblem2}) assembled only from $\LapP{n}$'s or $\LapN{n}$'s are singular
with simple zero eigenvalue. In all these cases, the corresponding eigenvector, i.e.,
the vector in null-space of such discretized Laplacian is
\[
  1_\ell = [1,1,\ldots,1]^\trp \in\mathbb{R}^\ell
\]
of suitable length $\ell$. Thus, the discretized right-hand side $H\in\mathbb{R}^{n\times q}$,
or $\mathcal{H}\in\mathbb{R}^{n\times q\times t}$ satisfies
\[
  \begin{split}
  1_{nq}\perp\vec(H) \qquad &\Longleftrightarrow \qquad 1_n^\trp H \,1_q = 0, \\
  1_{nqt}\perp\vec(\mathcal{H}) \qquad &\Longleftrightarrow \qquad (1_n^\trp,1_q^\trp,1_t^\trp\,|\,\mathcal{H}) = 0.
  \end{split}
\]
If it is not the case (in particular due to the rounding errors in the software that
provides us the right-hand side; but also during the computation, see Remark~\ref{rmk:orthogonality};
in principal also in the case, when the right-hand side is somehow measured), it is suitable
to project-out data in the direction of vector $1$, i.e., to do a modification of the form
\begin{equation}\label{eq:centering}
  \mbox{}\quad
  H \;\longmapsto\; H - 1_{n,q}\cdot\frac1{nq}\cdot\sum_{i=1}^n\sum_{j=1}^q h_{i,j}, \quad
  \mathcal{H} \;\longmapsto\; \mathcal{H} - 1_{n,q,t}\cdot\frac1{nqt}\cdot\sum_{i=1}^n\sum_{j=1}^q\sum_{k=1}^t h_{i,j,k}, \qquad
\end{equation}
i.e., to center the data.

\subsubsection{Dirichlet and Neumann BCs update}\label{sssec:RHS-DNBC}

If a Dirichlet or a Neumann BC is prescribed, then there is given 
particular value of electric potential or its derivative (electric field), constant on
the whole side or face. We demonstrate the right-hand side update for simplicity
only on the 1D case $L_n u=h$.

Consider the Dirichlet BC at the ``begining'' of 1D sample and let the prescribed
value of the electric potential be $u_{\mathsf{B}}$. Technically the prescription is done such that
\[
  \begin{array}{c}
  \mbox{\;} \\[-2pt]
  \left\lceil\!\begin{array}{c|ccccc}
    -1 & 2 & -1 \\
    & -1 & 2 & -1 \\
    & & \ddots & \ddots & \ddots \\
  \end{array}\!\right\rceil
  \end{array}
  \!\!
  \left\lceil\!\begin{array}{c}
    u_{\mathsf{B}} \\\hline
    u_1 \\ u_2 \\ \vdots
  \end{array}\!\right\rceil
  \!
  \begin{array}{c}
    \mbox{\;} \\ \mbox{\;} \\ \mbox{\;} \\ =
  \end{array}
  \!\!\!\!
  \begin{array}{c}
  \mbox{\;} \\[-2pt]
  \left\lceil\!\begin{array}{c}
    h_1 \\ h_2 \\ \vdots
  \end{array}\!\right\rceil
  \!\!\!
  \begin{array}{c}
    \mbox{\;} \\ \mbox{\;} \\ \text{,}
  \end{array}
  \end{array}
\]
which can be equivalently written as
\[
  \LapD{n}\,u = h + i_1u_{\mathsf{B}},
\]
where $i_\ell$ denotes the $\ell$th column of $I$.
(Similar update $i_n u_{\mathsf{E}}$ will be originated also from the value prescribed ``end''.)

Neumann BC at the ``begining'' of 1D sample prescribes the value of the electric field,
i.e., the derivative of electric potential over the boundary. Let it be $e_{\mathsf{B}}$. Technically
\[
  \begin{array}{c}
  \mbox{\;} \\[-2pt]
  \left\lceil\!\begin{array}{c|ccccc}
    -1 & 2 & -1 \\
    & -1 & 2 & -1 \\
    & & \ddots & \ddots & \ddots \\
  \end{array}\!\right\rceil
  \end{array}
  \!\!
  \left\lceil\!\begin{array}{c}
    u_1+e_{\mathsf{B}} \\\hline
    u_1 \\ u_2 \\ \vdots
  \end{array}\!\right\rceil
  \!
  \begin{array}{c}
    \mbox{\;} \\ \mbox{\;} \\ \mbox{\;} \\ =
  \end{array}
  \!\!\!\!
  \begin{array}{c}
  \mbox{\;} \\[-2pt]
  \left\lceil\!\begin{array}{c|ccccc}
    -1 & 1 & -1 \\
    & -1 & 2 & -1 \\
    & & \ddots & \ddots & \ddots \\
  \end{array}\!\right\rceil
  \end{array}
  \!\!
  \left\lceil\!\begin{array}{c}
    e_{\mathsf{B}} \\\hline
    u_1 \\ u_2 \\ \vdots
  \end{array}\!\right\rceil
  \!
  \begin{array}{c}
    \mbox{\;} \\ \mbox{\;} \\ \mbox{\;} \\ =
  \end{array}
  \!\!\!\!
  \begin{array}{c}
  \mbox{\;} \\[-2pt]
  \left\lceil\!\begin{array}{c}
    h_1 \\ h_2 \\ \vdots
  \end{array}\!\right\rceil
  \!\!\!
  \begin{array}{c}
    \mbox{\;} \\ \mbox{\;} \\ \text{,}
  \end{array}
  \end{array}
\]
which can be equivalently written as
\[
  \LapN{n}\,u = h + i_1e_{\mathsf{B}}.
\]
(Similar update $i_n e_{\mathsf{E}}$ will be originated also from the value prescribed ``end''.)

These two BCs can be easily combined in the 1D case, leading to $\LapDN{n}$ or $\LapND{n}$ matrix
and mixed updates $(i_1 u_{\mathsf{B}}+i_n e_{\mathsf{E}})$ or $(i_1 e_{\mathsf{B}}+i_n u_{\mathsf{E}})$
of the right-hand side.

\section{Preconditioned method of conjugate gradients}\label{sec:PCG}

Since the discretized Laplacians $L_n$ and in particular (\ref{eq:ourproblem2})
are either symmetric positive definite, or symmetric positive semi-definite, but
compatible, i.e., with the right-hand side orthogonal to the null-space of the operator
(see Section~\ref{sssec:RHS-nullspace}), we can use the method of conjugate gradients 
(CG), in particular its preconditioned variant (PCG) for solving such system;
see \cite{CGpaper}; see also
\cite{b:liesenstrakos12} for deeper understanding and wider context,
\cite{b:meurant06} for computational aspects,
\cite{b:malekstrakos15} for preconditioning of CG, and
\cite{b:meuranttichy24}, \cite{papez} for error estimation and stopping criterions.

Recall that the preconditioner is a linear operator. After discretization of the problem, preconditioner
represented by matrix $M$ is applied on both sides of $L_n u = h$ (in the 1D case), yielding 
the preconditioned system $ML_n u = Mh$. If $M$ is symmetric and positive definite, then there 
exists matrix $M^{\frac12}$ so that $M = M^{\frac12}M^{\frac12}$. The preconditioned system can 
be then further rearranged to  $(M^{\frac12}L_nM^{\frac12})(M^{-\frac12}u)=(M^{-\frac12}h)$ with 
symmetric matrix $M^{\frac12}L_nM^{\frac12}$ that is of the same definiteness as $L_n$. Consequently,
the standard CG can be applied to the latest system; however it approximates $M^{-\frac12}u$ vector.
After a few manipulations to express $u$, we get the usual PCG algorithm.

\subsection{PCG algorithm}

Since the PCG is a very standard algorithm, we explicitly mention only its variant
for our 3D problem; see Algorithm \ref{alg:PCG:3D}. To get its 2D variant it is sufficient
to replace all the three-way tensors by matrices --- formally replace all the calligraphic
letters $\mathcal{H}$, $\mathcal{U}_s$, $\mathcal{W}_s$, $\mathcal{R}_s$,
$\mathcal{Z}_s$, and $\mathcal{P}_s$, by its uppercase italic counterparts.

\begin{algorithm}[htb!]
\caption{Preconditioned CG for discretized 3D Poisson's equation} \label{alg:PCG:3D}
\begin{algorithmic}[1]
\STATE\label{line01} {\bf Input} $\FDL$, $\mathcal{H}$, $\mathfrak{M}$, $\mathcal{U}_0$
  \hfill\COMMENT{Laplacian, right-hand side, preconditioner, initial guess}
\STATE\label{line02} $s \leftarrow 0$
\STATE\label{line03} $\mathcal{W}_0 \leftarrow -\FDL(\mathcal{U}_0)$
  \hfill\COMMENT{application of discretized minus Laplacian}
\STATE\label{line04} $\mathcal{R}_0 \leftarrow \mathcal{H} - \mathcal{W}_0$
  \hfill\COMMENT{initial residuum}
\STATE\label{line05} $\mathcal{Z}_0 \leftarrow \mathfrak{M}(\mathcal{R}_0)$
  \hfill\COMMENT{application of preconditioner}
\STATE\label{line06} $\rho_0 \leftarrow \langle\mathcal{R}_0,\mathcal{Z}_0\rangle$
\STATE\label{line07} $\mathcal{P}_0 \leftarrow \mathcal{Z}_0$
\STATE\label{line08} {\bf Repeat}
\STATE\label{line09} \qquad $s \leftarrow s+1$
\STATE\label{line10} \qquad $\mathcal{W}_{s} \leftarrow -\FDL(\mathcal{P}_{s-1})$
  \hfill\COMMENT{application of discretized minus Laplacian}
\STATE\label{line11} \qquad $\alpha_{s} \leftarrow \rho_{s-1}/\langle\mathcal{W}_{s},\mathcal{P}_{s-1}\rangle$
\STATE\label{line12} \qquad $\mathcal{U}_{s} \leftarrow \mathcal{U}_{s-1} + \alpha_{s} \mathcal{P}_{s-1}$
  \hfill\COMMENT{update of solution approximation}
\STATE\label{line13} \qquad $\mathcal{R}_{s} \leftarrow \mathcal{R}_{s-1} - \alpha_{s} \mathcal{W}_{s}$
  \hfill\COMMENT{update of residuum}
\STATE\label{line14} \qquad $\mathcal{Z}_{s} \leftarrow \mathfrak{M}(\mathcal{R}_{s})$
  \hfill\COMMENT{application of preconditioner}
\STATE\label{line15} \qquad $\rho_{s} \leftarrow \langle\mathcal{R}_{s},\mathcal{Z}_{s}\rangle$
\STATE\label{line16} \qquad $\beta_{s} \leftarrow \rho_{s}/\rho_{s-1}$
\STATE\label{line17} \qquad $\mathcal{P}_{s} \leftarrow \mathcal{Z}_{s} + \beta_{s} \mathcal{P}_{s-1}$
\STATE\label{line18} {\bf Until} stopping criterion match
\STATE\label{line19} {\bf Output} $\mathcal{U}_{\mathsf{PCG}} \leftarrow \mathcal{U}_{s}$
  \hfill\COMMENT{PCG approximation of $\mathcal{U}$}
\end{algorithmic}
\end{algorithm}


\begin{remark}[Orthogonality to null-space]\label{rmk:orthogonality}
Consider the discretized Laplacian which is singular (in 1D case $\LapP{n}$,
see Lemma~\ref{lem:spectrumP}, or $\LapN{n}$, see Lemma \ref{lem:spectrumN}).
Its null-space is one-dimensional, spaned by $1_n=[1,1,\dots,1]^\trp$ in 1D,
$1_{n,q}$ in 2D, and $1_{n,q,t}$ in 3D case.

If the right-hand side $\mathcal{H}$ is orthogonal to the null-space (see
Section~\ref{sssec:RHS-nullspace}), and also the initial guess $\mathcal{U}_0$
(notice that the initial guess is often zero, $\mathcal{U}_0=0_{n,q,t}$), then
the whole CG (PCG with $\mathfrak{M}$ being identity) is orthogonal to this
null-space. In particular vectors, $\mathcal{U}_s$, $\mathcal{W}_s$,
$\mathcal{R}_s$, and $\mathcal{P}_s$ are orthogonal to $1_{n,q,t}$.
In the finite precision arithmetic, a nonzero projection into the null-space
may be generated due to rounding errors. Therefore, it is suitable to center
the data; see (\ref{eq:centering}). It seems to be sufficient to center only
one of the four quantities, once per iteration.

In the case of general PCG, the preconditioner $\mathfrak{M}$ may generate
nonzero projection of $\mathcal{Z}_s$ into the null-space and thus affect
the whole algorithm. Therefore, it is suitable to center $\mathcal{Z}_s$.
\end{remark}


\begin{remark}[Stopping criterion]\label{rmk:stoppincrit}
Since this work is rather a case study about behavior the PCG with different 
preconditioners when applied on our problem, we do not elaborate the stopping 
criterion. Actually, we stop simply at maximal number of iterations. This number 
is determined experimentally for each problem and preconditioner, so that the 
true residual decreases sufficiently. This choice seems to be, however, applicable
even in our final practical C/C++ implementation.

The best option for stopping is robust estimation of $A$-norm of error; see \cite{b:meuranttichy24}
and in particular paper \cite{papez}.
\end{remark}

\subsection{Computational cost of PCG}

The PCG algorithm consists of several typical subroutines. Some of them needs to
be explored in details. We start from the simplest:
\begin{itemize}\itemsep 0pt
\item {\bf Saxpy type subroutine} (lines \ref{line04}, \ref{line12}, \ref{line13}, and \ref{line17})
that consists of simple entry-wise operations with matrices or three-way tensors.
The cost of one saxpy is about $O(2nqt)$ of elementary operations; in the 2D case $t=1$.
\item {\bf Inner product} (lines \ref{line06}, \ref{line11}, and \ref{line15}) that can be for matrices
and three-way tensors defined by employing the vectorization, e.g.,
\[
  \begin{split}
  \langle R,Z\rangle = \langle\vec(R),\vec(Z)\rangle
  &= \sum_{i=1}^n\sum_{j=1}^q {r}_{i,j}\cdot{z}_{i,j},
  \\
  \langle\mathcal{R},\mathcal{Z}\rangle = \langle\vec(\mathcal{R}),\vec(\mathcal{Z})\rangle
  &= \sum_{i=1}^n\sum_{j=1}^q\sum_{k=1}^t {r}_{i,j,k}\cdot{z}_{i,j,k}.
  \end{split}
\]
The cost of one inner product is about $O(2nqt)$ of elementary operations.
\item {\bf Application of minus Laplacian} (lines \ref{line03} and \ref{line10}) that needs to be broken
down to operations with the small 1D Laplacian matrices $L_\ell$, e.g.,
\[
  \mbox{}\qquad\qquad
  W = -\FDL(P) = L_n P + P L_q^\trp,
  \quad
  \mathcal{W} = -\FDL(\mathcal{P}) = L_n\times_1\mathcal{P} + L_q\times_2\mathcal{P} + L_t\times_3\mathcal{P},
\]
with entries
\[
  \begin{split}
  w_{\zeta,\xi} &= \sum_{i=1}^n (L_n)_{\zeta,i}\cdot p_{i,\xi} + \sum_{j=1}^q (L_q)_{\xi,j}\cdot p_{\zeta,j},
  \\
  w_{\zeta,\xi,\nu} &= \sum_{i=1}^n (L_n)_{\zeta,i}\cdot p_{i,\xi,\nu} + \sum_{j=1}^q (L_q)_{\xi,j}\cdot p_{\zeta,j,\nu} + \sum_{k=1}^t (L_t)_{\nu,k}\cdot p_{\zeta,\xi,k}.
  \end{split}
\]
The cost of one application of a general operator with this structure is about $O(2nq(n+q))$ and $O(2nqt(n+q+t))$,
in 2D and 3D case, respectively. Here, however, the sparsity of $L_\ell$ matrices can (and needs to) be utilized
in the practical implementation --- since they have only three nonzero entries per row, the cost can be reduced
to about $O(12nqt)$ and $O(18nqt)$, respectively.
\item {\bf Application of preconditioner} (lines \ref{line05} and \ref{line14}) depends
on the particular choice of preconditioner. Clearly, it is useful to look for such one,
which do not increase the computational cost significantly.
\end{itemize}

The overall cost of our PCG implementation is in Table \ref{tab:PCGcost}. The cost of
preconditioning is not taken into the account. In the case of singular discretized Laplacian,
each data centering (see Section~\ref{sssec:RHS-nullspace} and Remark~\ref{rmk:orthogonality})
increases the cost by $3nqt$ elementary operations --- typically once per iteration. Similarly,
the evaluation of stopping criterion (see Remark~\ref{rmk:stoppincrit}) may costs some further operations.

\begin{table}[thb!]
\caption{Cost (measured by number of elementary operations) of matrix (2D) and tensor (3D)
variants of PCG with Kronecker-structured operator with sparse matrices (three nonzeros per row);
see Algorithm \ref{alg:PCG:3D}. Preconditioner is not included.}\label{tab:PCGcost}
\begin{tabular}{lllll}
                   & \multicolumn{2}{l}{Initialization (lines \ref{line02}--\ref{line07})}
                   & \multicolumn{2}{l}{Iteration (lines \ref{line09}--\ref{line17})} \\\hline
App. of Laplacian  & Full            & Sparse  & Full            & Sparse  \\\hline
Cost of 2D variant & $2nq(n+q+2)$    & $16nq$  & $2nq(n+q+5)$    & $22nq$  \\
Cost of 3D variant & $2nqt(n+q+t+2)$ & $22nqt$ & $2nqt(n+q+t+5)$ & $28nqt$ \\
\end{tabular}
\end{table}

\subsection{Preconditioner}

Typically we are looking for such preconditioner, that somehow resembles the tensor (or matrix)
structure of the algorithm, which does not increase the computational cost, and which, of course,
speed up the convergence. Notice that relevant is the computational cost per iteration, not the
cost of the preparatory work, that is done once, at the beginning of PCG.
We have experimented with only few families of preconditioners:
\begin{itemize}\itemsep 0pt
\item Several steps of stationary iterative method --- because of symmetry of $L_\ell$s,
we decided to use Jacobi-like method (standard and weighted Jacobi).
\item Approximate inverse of $(-\FDL)$ of low-Kronecker-rank (Kronecker-rank
denotes the number of summands when the matrix is expressed as in (\ref{eq:ourproblem2});
notice that $(-\FDL)$ is always of low-Kronecker-rank, equal to two or three,
in the 2D and 3D case, respectively).
\item Moore--Penrose pseudoinverse of $(-\FDL)$.
\end{itemize}

Implementation of a preconditioner into PCG consists of two subroutines:
The first one prepares all the parts of the preconditioner that can be
computed only once, in beforehand, while initialization of PCG.
The other one is the application of the precomputed preconditioner on the
updated residual (see lines \ref{line05} and \ref{line14} of Algorithm \ref{alg:PCG:3D}).

\subsubsection{Preconditioning by several steps of Jacobi-like method}\label{sssec:SSJLM}

Consider parametrized splitting of matrix $A\in\mathbb{R}^{n\times n}$ in the form
\[
  A = D_\omega(A) + O_\omega(A),
  \qquad \text{where} \qquad
  D_\omega(A) = D\cdot\omega, \quad
  O_\omega(A) = O + D\cdot(1-\omega),
\]
and where $D=\diag(a_{1,1},a_{2,2},\ldots,a_{n,n})\in\mathbb{R}^{n\times n}$
contains the diagonal entries of $A$, $O = A - D$ the off-diagonal entries, and
$\omega\in[1,\infty)$ is the real weighting parameter. The associated stationary
method applied on $Ax=b$ iterates as
\begin{equation}\label{eq:jacobimethod}
  x_{j} = (D_\omega(A))^{-1}(b-O_\omega(A)x_{j-1}).
\end{equation}
If $\omega=1$, the splitting and the associated method is called Jacobi,
if $\omega>1$, it is called weighted Jacobi.

Matrix $D_\omega(A)$ is diagonal, thus its inverse is easy to compute --- just
invert all the diagonal entries. If we handle its diagonal as a vector
\[
  \hat{d}_\omega = \left[\begin{array}{c}
  a_{1,1} \\
  a_{2,2} \\
  \vdots  \\
  a_{n,n}
  \end{array}\right]\cdot\omega \in\mathbb{R}^n \qquad \text{and denote} \qquad
  f_{j} = b-O_\omega(A)x_{j-1}\in\mathbb{R}^n,
\]
then we can (\ref{eq:jacobimethod}) rewrite as
\begin{equation}\label{eq:jacobihadamard}
  x_{j} = (\hat{d}_\omega)^{\odot-1}\odot f_j.
\end{equation}
where ${}^{\odot-1}$ denotes the {\em entry-wise inverse} also called {\em Hadamard inverse}
of the vector, and $\odot$ the {\em entry-wise} a.k.a. {\em Hadamard product}.

For Kronecker-structured matrix such as in (\ref{eq:ourproblem2}), e.g.,
\[
  \LapDD = (I_q\otimes L_n)+(L_q\otimes I_n)
\]
in the 2D case, it is easy to see (using property (\ref{eq:kronPsum})) that the
splitting takes the same form, i.e.,
\begin{equation}\label{eq:DOkron}
  \begin{split}
  D_\omega(\LapDD) &= (I_q\otimes D_\omega(L_n))+(D_\omega(L_q)\otimes I_n), \\
  O_\omega(\LapDD) &= (I_q\otimes O_\omega(L_n))+(O_\omega(L_q)\otimes I_n).
  \end{split}
\end{equation}
Iteration (\ref{eq:jacobimethod}) requires inverse of diagonal matrix $D_\omega(\LapDD)$,
which can be again realized as in (\ref{eq:jacobihadamard}), by Hadamard inverse of
\begin{equation}\label{eq:HadInv}
  \hat{D}_\omega = \left[\begin{array}{cccc}
    (L_n)_{1,1}+(L_q)_{1,1} & (L_n)_{1,1}+(L_q)_{2,2} & \cdots & (L_n)_{1,1}+(L_q)_{q,q} \\
    (L_n)_{2,2}+(L_q)_{1,1} & (L_n)_{2,2}+(L_q)_{2,2} & \cdots & (L_n)_{2,2}+(L_q)_{q,q} \\
    \vdots & \vdots & \ddots & \vdots \\
    (L_n)_{n,n}+(L_q)_{1,1} & (L_n)_{n,n}+(L_q)_{2,2} & \cdots & (L_n)_{n,n}+(L_q)_{q,q}
  \end{array}\right]\cdot\omega
\end{equation}
and then by Hadamard product.
To summarize, the only parts that needs to be precomputed are: Matrices
$O_\omega(L_n)$, $O_\omega(L_q)$, and $(\hat{D}_\omega)^{\odot-1}$ in the 2D case;
matrices $O_\omega(L_n)$, $O_\omega(L_q)$, $O_\omega(L_t)$, and
tensor $(\hat{\mathcal{D}}_\omega)^{\odot-1}$ in the 3D case.

For preconditioning CG we use $p$ steps of this stationary method. It is applied
on the problem with the updated residual on the right-hand side (see lines \ref{line05}
and \ref{line14} of Algorithm \ref{alg:PCG:3D}) and with initial guess $x_{0}=0$.
See Algorithm~\ref{alg:SSJLM:appl:2D} for implementation for the 2D case and
Algorithm~\ref{alg:SSJLM:appl:3D} for the 3D case.

\begin{algorithm}[htb!]
\caption{Application of preconditioner $Z_s \leftarrow \mathfrak{M}_{\text{Jacobi}(p,\omega)}(R_s)$, 2D case} \label{alg:SSJLM:appl:2D}
\begin{algorithmic}[1]
\STATE\label{P1-2D-line01} {\bf Precomputed} $O_\omega(L_n)$, $O_\omega(L_q)$, $(\hat{D}_\omega)^{\odot-1}$
\STATE\label{P1-2D-line02} {\bf Input} $R_s$
  \hfill\COMMENT{updated residual}
\STATE\label{P1-2D-line03} $X_0 = 0_{n,q}$
\STATE\label{P1-2D-line04} {\bf For} $j \leftarrow 1,2,\ldots,p$
\STATE\label{P1-2D-line05} \qquad $F_{j} \leftarrow R_s - (O_\omega(L_n)X_{j-1} + X_{j-1}(O_\omega(L_q))^\trp)$
\STATE\label{P1-2D-line06} \qquad $X_{j} \leftarrow (\hat{D}_\omega)^{\odot-1} \odot F_{j}$
\STATE\label{P1-2D-line07} {\bf End for}
\STATE\label{P1-2D-line08} {\bf Output} $Z_s \leftarrow X_{p}$
  \hfill\COMMENT{preconditioned residual}
\end{algorithmic}
\end{algorithm}

\begin{algorithm}[htb!]
\caption{Application of preconditioner $\mathcal{Z}_s \leftarrow \mathfrak{M}_{\text{Jacobi}(p,\omega)}(\mathcal{R}_s)$, 3D case} \label{alg:SSJLM:appl:3D}
\begin{algorithmic}[1]
\STATE\label{P1-3D-line01} {\bf Precomputed} $O_\omega(L_n)$, $O_\omega(L_q)$, $O_\omega(L_t)$, $(\hat{\mathcal{D}}_\omega)^{\odot-1}$
\STATE\label{P1-3D-line02} {\bf Input} $\mathcal{R}_s$
  \hfill\COMMENT{updated residual}
\STATE\label{P1-3D-line03} $\mathcal{X}_0 = 0_{n,q,t}$
\STATE\label{P1-3D-line04} {\bf For} $j \leftarrow 1,2,\ldots,p$
\STATE\label{P1-3D-line05} \qquad $\mathcal{F}_{j} \leftarrow \mathcal{R}_s - (O_\omega(L_n)\times_1\mathcal{X}_{j-1} + O_\omega(L_q)\times_2\mathcal{X}_{j-1} + O_\omega(L_t)\times_3\mathcal{X}_{j-1})$
\STATE\label{P1-3D-line06} \qquad $\mathcal{X}_{j} \leftarrow (\hat{\mathcal{D}}_\omega)^{\odot-1} \odot \mathcal{F}_{j}$
\STATE\label{P1-3D-line07} {\bf End for}
\STATE\label{P1-3D-line08} {\bf Output} $\mathcal{Z}_s \leftarrow \mathcal{X}_{p}$
  \hfill\COMMENT{preconditioned residual}
\end{algorithmic}
\end{algorithm}

\subsubsection{Preconditioning by a low-Kronecker-rank approximation of the inverse}\label{sssec:AILKR}

Recall that (\ref{eq:ourproblem}) can be rewritten as standard system (\ref{eq:ourproblem2})
with matrix $((I_q\otimes L_n)+(L_q\otimes I_n))$. We would like to approximate its inverse,
or generalized inverse (see \cite{geninv}) since the original matrix may be singular, by a low
Kronecker-rank matrix. In particular
\[
  \LapDD^\dagger
  = \Big((I_q\otimes L_n)+(L_q\otimes I_n)\Big)^{\dagger}
  \approx \sum_{\rho=1}^{r_\kro}
  (M_{q,\rho} \otimes M_{n,\rho}),
  \qquad \text{where} \qquad
  M_{\ell,\rho} \in\mathbb{R}^{\ell\times \ell}
\]
and where
${}^\dagger$ denotes the Moore--Penrose pseudoinverse (see \cite{geninv}, \cite[Section 5.2.2]{GolubVanLoan}) and
$r_\kro$ is low (i.e., comparable with the Kronecker-rank of the original
matrix, i.e., in the 2D case $r_\kro\lesssim2$); for related discussion see \cite{KronPrec1}, \cite{KronPrec2}.

Recall the eigenvalue decompositions of $L_\ell$s (see Lemmas \ref{lem:spectrumP}--\ref{lem:spectrumM}).
Notice that all these matrices are positive semi-definite, i.e., their eigenvalues are nonnegative.
If we normalize all the eigenvectors of $L_\ell$, then the decomposition (\ref{eq:eigL2D1}) represent
the singular value decomposition (SVD) of $L_\ell$,
\begin{equation}\label{eq:eig-svd-decomp}
  L_\ell = V_\ell\Lambda_\ell V_\ell^{\trp}, \qquad V_\ell^{-1} = V_\ell^{\trp}.
\end{equation}
Now, since the right-hand side $H\in\mathbb{R}^{n\times q}$ is in general dense matrix without any particular
structure explicitly exploitable for compression, we need to be able to store several such matrices
(i.e., of $\sim nq$ entries) to proceed the PCG algorithm. Thus, if $n\sim q$, then also $n^2\sim q^2\sim nq$,
and then we are also able to store the full eigendecompositions/SVDs of all $L_\ell$s.
Moreover, we know these decompositions analytically.\footnote{However it seems, that the eigendecompositions
computed in {\sc Matlab} by {\tt eig} function works better than the decompositions assembled entry-wise
by using Lemmas \ref{lem:spectrumP}--\ref{lem:spectrumM}. Better means the yielding a smaller
Frobenius norm $\|L_\ell - V_\ell\Lambda_\ell V_\ell^{\trp}\|_\fro$.}

Consequently, we have all the important components for the eigendecomposition/SVD of
the big matrix $\LapDD$ (\ref{eq:eigL2D2}),
\[
  \Big((I_q\otimes L_n)+(L_q\otimes I_n)\Big)
  = (V_q\otimes V_n)\,\Big((I_q\otimes \Lambda_n)+(\Lambda_q\otimes I_n)\Big)\,(V_q\otimes V_n)^{\trp},
\]
and also for its pseudoinverse
\begin{equation}\label{eq:p-inv-full}
  \Big((I_q\otimes L_n)+(L_q\otimes I_n)\Big)^\dagger
  = (V_q\otimes V_n)\,\Big((I_q\otimes \Lambda_n)+(\Lambda_q\otimes I_n)\Big)^\dagger\,(V_q\otimes V_n)^{\trp}.
\end{equation}
The only thing that remains is to calculate low-Kronecker-rank approximation
of the matrix in the middle. This matrix is however diagonal
\[
  \Big((I_q\otimes \Lambda_n)+(\Lambda_q\otimes I_n)\Big) = \diag(\,\cdots\;,\;\lambda_i(L_n)+\lambda_j(L_q)\;,\;\cdots\,),
\]
thus its pseudoinverse takes form
\[
  \Big((I_q\otimes \Lambda_n)+(\Lambda_q\otimes I_n)\Big)^\dagger = \diag(\,\cdots\;,\;g_{i,j}\;,\;\cdots\,),
  \quad \text{where} \quad
  g_{i,j} = \left\{\begin{array}{c}
    \frac{1}{\lambda_i(L_n)+\lambda_j(L_q)} \\ 0
  \end{array}\right.,
\]
i.e., the nonzero diagonal entries are inverted, the zeros are left as they are.
Notice that in practical implementation we consider all the diagonal entries with
modulus smaller than $10^{-13}$ to be zero, thus we replace these entries by zeros
in the pseudoinverse.

Now we use similar approach as in previous Section \ref{sssec:SSJLM}. The inverse
of diagonal Kronecker-structured matrix $D_\omega$ (\ref{eq:DOkron}) can be written
as the Hadamard inverse of (\ref{eq:HadInv}). Thus,
\begin{equation}\label{eq:had-p-inv}
  \hat{G} = \left[\begin{array}{ccc}
  g_{1,1} & \cdots & g_{1,q} \\
  \vdots & \ddots & \vdots \\
  g_{n,1} & \cdots & g_{n,q}
  \end{array}\right]
  = \left[\begin{array}{ccc}
  \lambda_1(L_n)+\lambda_1(L_q) & \cdots & \lambda_1(L_n)+\lambda_q(L_q) \\
  \vdots & \ddots & \vdots \\
  \lambda_n(L_n)+\lambda_1(L_q) & \cdots & \lambda_n(L_n)+\lambda_q(L_q)
  \end{array}\right]^{\odot\dagger},
\end{equation}
where ${}^{\odot\dagger}$ is the Hadamard (entry-wise) pseudoinverse.

The remaining steps are rather straightforward. To get the best Kronecker-rank $r_\kro$
approximation of pseudoinverse of the orignal diagonal matrix, we need to compute
the best rank $r_\kro$ (in the standard sense) approximation of $\hat{G}$.
This can be obtained by using the SVD of $\hat{G}$, and employing the Eckart--Young--Mirsky
theorem \cite{EckartYoung}, \cite{Mirsky}. Let
\[
  \hat{G}=X\Sigma Y^\trp = \sum_{\rho=1}^{\min(n,d)} x_\rho\sigma_\rho y_\rho^\trp
\]
be the SVD and its the dyadic expansion. The first $r_\kro$ dyads, i.e., outer
products ${x_\rho}{\sigma_\rho}{y_\rho}^\trp$ corresponding to largest singular
values can be easily reshaped to Kronecker products of diagonal matrices, e.g.,
as $\diag({x_\rho}{\sigma_\rho})\otimes\diag({y_\rho})$. This gives us
\[
  \begin{split}
  \LapDD^\dagger =
  \Big((I_q\otimes L_n)+(L_q\otimes I_n)\Big)^\dagger
  &\approx (V_q\otimes V_n)\,\Big(\sum_{\rho=1}^{r_\kro} \diag({x_\rho}{\sigma_\rho})\otimes\diag({y_\rho})\Big)\,(V_q\otimes V_n)^{\trp} \\
  &= \sum_{\rho=1}^{r_\kro}
    \Big(\underbrace{V_q\,\diag({x_\rho}{\sigma_\rho})V_q^\trp}_{\displaystyle M_{q,\rho}}\Big) \otimes
    \Big(\underbrace{V_n\,\diag({y_\rho})V_n^\trp}_{\displaystyle M_{n,\rho}}\Big),
  \end{split}
\]
the wanted approximation.

We do not provide the algorithmic description of implementation of this preconditioner here,
for several reasons. The {\em first} of them may be more-or-less obvious: If we are able to
compute, store, and apply the whole matrix $\hat{G}\in\mathbb{R}^{n\times q}$ --- and we are,
because we are able to compute all the matrices $\Lambda_\ell$, to store several matrices
of $\sim nq$ entries, and its application is done simply by the Hadamard product, as in
Section \ref{sssec:SSJLM}  --- we do not need to approximate this matrix by anything, nor by
low-rank object. The low-Kronecker-rank approach seems to be superfluous, it does not work
as good as full Kronecker-rank, with about the same cost; see the next Section \ref{sssec:MPPI}.
The {\em second} related reason: Initialization of this preconditioner is almost identical
to the next one (see Algorithm \ref{alg:MPPI:init:2D}); here we need to extra precompute the SVD
of $\hat{G}$ and prepare matrices $M_{n,\rho}$, $M_{q,\rho}$. Application of this preconditioner
(see lines \ref{line05} and \ref{line14} of Algorithm \ref{alg:PCG:3D}) is almost identical to
the application of Laplacian, since both have the same structure, in particular
\[
  Z_s =
  \mathfrak{M}_{\text{low-rank}(r_\kro)}(R_s) =
  \sum_{\rho=1}^{r_\kro} M_{n,\rho} R_s M_{q,\rho}^\trp.
\]
The {\em third} reason is less obvious: There is no simple generalization of this approach to
the 3D case. The low-rank approximation of $\hat{\mathcal{G}}\in\mathbb{R}^{n\times q\times t}$
needs to be done by some kind of high-order SVD (HOSVD), typically the Tucker decomposition
(see \cite{Tucker1}, \cite{Tucker2}, \cite{Tucker3}) or the CP (CanDeComp \cite{Candecomp} $+$
ParaFac \cite{Parafac}) decomposition; see also \cite{KolBad09}. Tucker decomposition provides
optimal approximation, however, the so-called Tucker-core --- an analogy of the $\Sigma$ matrix
from the SVD --- is not diagonal, but in general full. Thus the total number of the Kronecker
summands grows rapidly (cubic in comparison to the SVD in the 2D case). On the other hand,
CP decomposition provides diagonal core tensor and triadic (in general polyadic) decomposition
similarly as the SVD, however, the so-called polyadic rank which is involved here, is not well
posed due to the lack of orthogonality among individual triads.
The {\em fourth}, the last but least reason is the positive semi-definiteness of
$\mathfrak{M}_{\text{low-rank}(r_\kro)}$. Its eigenvalues are entries of the partial
sum $\sum_{\rho=1}^{r_\kro} {x_\rho}{\sigma_\rho}{y_\rho}^\trp$. Since vectors $x_\rho$
are mutually orthonormal, and $y_\rho$ as well, there is no way (except for
$r_\kro = 1$ and $r_\kro = \min(n,q)$) to guarantee these entries to be nonnegative.

\subsubsection{Preconditioning by Moore--Penrose pseudoinverse}\label{sssec:MPPI}

As already indicated in the previous section, to get all necessary components
the Moore--Penrose pseudoinverse of $(-\FDL)$, we only need to compute the
eigendecompositions/SVDs of individual small matrices $L_\ell$ (\ref{eq:eig-svd-decomp}),
which is manageable, and then calculate the Hadamard pseudoinverse of the matrix/tensor
of all sums of all eigenvalues of different $L_\ell$s (\ref{eq:had-p-inv}).
Recall that we consider all the sums with modulus smaller than $10^{-13}$
to be zero in practical implementation. This is all we need to initialize this
preconditioner; see also Algorithms \ref{alg:MPPI:init:2D} and \ref{alg:MPPI:init:3D}.

\begin{algorithm}[htb!]
\caption{Initialization of preconditioner $\mathfrak{M}_{\text{p-inv}}$, 2D case} \label{alg:MPPI:init:2D}
\begin{algorithmic}[1]
\STATE\label{P3-2D-init-line01} {\bf Input} $L_n$, $L_q$
  \hfill\COMMENT{1D Laplacians}
\STATE\label{P3-2D-init-line02} $[V_n,\lambda_1(L_n),\ldots,\lambda_n(L_n)] \leftarrow \mathtt{eig}(L_n)$
\STATE\label{P3-2D-init-line03} $[V_{q\,},\lambda_1(L_{q\,}),\ldots,\lambda_{q\,}(L_{q\,})] \leftarrow \mathtt{eig}(L_q)$
\STATE\label{P3-2D-init-line04} {\bf For} $i \leftarrow 1,2,\ldots,n$ {\bf and} $j \leftarrow 1,2,\ldots,q$
\STATE\label{P3-2D-init-line05} \qquad $\hat{G}_{i,j} \leftarrow (\lambda_i(L_n) + \lambda_j(L_q))^\dagger$
  \hfill\COMMENT{threshold for nonzeros: $10^{-13}$}
\STATE\label{P3-2D-init-line06} {\bf End for}
\STATE\label{P3-2D-init-line07} {\bf Output} $V_n$, $V_q$, $\hat{G}$
  \hfill\COMMENT{eigenvectors, pseudoinverse of sums of eigenvalues}
\end{algorithmic}
\end{algorithm}

\begin{algorithm}[htb!]
\caption{Initialization of preconditioner $\mathfrak{M}_{\text{p-inv}}$, 3D case} \label{alg:MPPI:init:3D}
\begin{algorithmic}[1]
\STATE\label{P3-3D-init-line01} {\bf Input} $L_n$, $L_q$, $L_t$
  \hfill\COMMENT{1D Laplacians}
\STATE\label{P3-3D-init-line02} $[V_n,\lambda_1(L_n),\ldots,\lambda_n(L_n)] \leftarrow \mathtt{eig}(L_n)$
\STATE\label{P3-3D-init-line03} $[V_{q\,},\lambda_1(L_{q\,}),\ldots,\lambda_{q\,}(L_{q\,})] \leftarrow \mathtt{eig}(L_q)$
\STATE\label{P3-3D-init-line04} $[V_{t\,\,},\lambda_1(L_{t\,\,}),\ldots,\lambda_{t\,\,}(L_{t\,\,})] \leftarrow \mathtt{eig}(L_t)$
\STATE\label{P3-3D-init-line05} {\bf For} $i \leftarrow 1,2,\ldots,n$ {\bf and} $j \leftarrow 1,2,\ldots,q$ {\bf and} $k \leftarrow 1,2,\ldots,t$
\STATE\label{P3-3D-init-line06} \qquad $\hat{\mathcal{G}}_{i,j,k} \leftarrow (\lambda_i(L_n) + \lambda_j(L_q) + \lambda_k(L_t))^\dagger$
  \hfill\COMMENT{threshold for nonzeros: $10^{-13}$}
\STATE\label{P3-3D-init-line07} {\bf End for}
\STATE\label{P3-3D-init-line08} {\bf Output} $V_n$, $V_q$, $V_t$, $\hat{\mathcal{G}}$
  \hfill\COMMENT{eigenvectors, pseudoinverse of sums of eigenvalues}
\end{algorithmic}
\end{algorithm}

Application of this preconditioner can be easily understood from the structure of
pseudoinverse (\ref{eq:p-inv-full}). We take the updated residual, transform it
into eigenvector-bases, then divide by sums eigenvalues (multiply by their pseudoinverses),
and transform back; see also Algorithms \ref{alg:MPPI:appl:2D} and \ref{alg:MPPI:appl:3D}.

\begin{algorithm}[htb!]
\caption{Application of preconditioner $Z_s \leftarrow \mathfrak{M}_{\text{p-inv}}(R_s)$, 2D case} \label{alg:MPPI:appl:2D}
\begin{algorithmic}[1]
\STATE\label{P3-2D-appl-line01} {\bf Precomputed} $V_n$, $V_q$, $\hat{G}$
\STATE\label{P3-2D-appl-line02} {\bf Input} $R_s$
  \hfill\COMMENT{updated residual}
\STATE\label{P3-2D-appl-line03} $F_1 \leftarrow V_n^\trp R_s V_q$
\STATE\label{P3-2D-appl-line04} $F_2 \leftarrow F_1 \odot \hat{G}$
\STATE\label{P3-2D-appl-line05} $F_3 \leftarrow V_n F_2 V_q^\trp$
\STATE\label{P3-2D-appl-line06} {\bf Output} $Z_s\leftarrow F_3$
  \hfill\COMMENT{preconditioned residual}
\end{algorithmic}
\end{algorithm}

\begin{algorithm}[htb!]
\caption{Application of preconditioner $\mathcal{Z}_s \leftarrow \mathfrak{M}_{\text{p-inv}}(\mathcal{R}_s)$, 3D case} \label{alg:MPPI:appl:3D}
\begin{algorithmic}[1]
\STATE\label{P3-3D-appl-line01} {\bf Precomputed} $V_n$, $V_q$, $V_t$, $\hat{\mathcal{G}}$
\STATE\label{P3-3D-appl-line02} {\bf Input} $\mathcal{R}_s$
  \hfill\COMMENT{updated residual}
\STATE\label{P3-3D-appl-line03} $\mathcal{F}_1 \leftarrow (V_n^\trp,V_q^\trp,V_t^\trp\,|\,\mathcal{R}_s)$
\STATE\label{P3-3D-appl-line04} $\mathcal{F}_2 \leftarrow \mathcal{F}_1 \odot \hat{\mathcal{G}}$
\STATE\label{P3-3D-appl-line05} $\mathcal{F}_3 \leftarrow (V_n,V_q,V_t\,|\,\mathcal{F}_2)$
\STATE\label{P3-3D-appl-line06} {\bf Output} $\mathcal{Z}_s\leftarrow \mathcal{F}_3$
  \hfill\COMMENT{preconditioned residual}
\end{algorithmic}
\end{algorithm}

Table \ref{tab:PInvPcost} summarizes computational cost of this preconditioner while
complementing Table \ref{tab:PCGcost}. We see, that the cost of the pseudoinverse
preconditioner is comparable with the `full' (i.e., not the sparse variant) PCG itself.
Cost of the eigendecomposition of symmetric matrix of order $n$ is about $n^3$, moreover,
we deal with matrices (\ref{eq:lap1d}) that are either banded --- even tridiagonal
(with the exception of the case with periodic BC), or Toeplitz (with the exception of the
cases with Neumann BC), so the cost can be further reduced; see \cite[Chapter 8]{GolubVanLoan},
\cite{Parlett}, \cite{Saad}. If we employ the analytic
formulas given in Lemmas~\ref{lem:spectrumP}--\ref{lem:spectrumM}, we can reduce it
to about $n^2$.

\begin{table}[thb!]
\caption{Cost of the matrix (2D) and tensor (3D) variants of the Moore--Penrose pseudoinverse
preconditioner for PCG. The {\tt eig}\,s may be assembled directly,
or computed by an dedicated solver.}\label{tab:PInvPcost}
\begin{tabular}{lllll}
                   & \multicolumn{3}{l}{Initialization (Algs. \ref{alg:MPPI:init:2D}, \ref{alg:MPPI:init:3D})}
                   & \multicolumn{1}{l}{Application (Algs. \ref{alg:MPPI:appl:2D}, \ref{alg:MPPI:appl:3D})} \\\hline
Cost of 2D variant & Cost of $2$ {\tt eig}\,s & \!\!$+$\!\! & $2nq$  & $4nq(n+q+\frac14)$     \\
Cost of 3D variant & Cost of $3$ {\tt eig}\,s & \!\!$+$\!\! & $3nqt$ & $4nqt(n+q+t+\frac14)$  \\
\end{tabular}
\end{table}


\begin{remark}\label{rmk:fourier}
First it is worth to see once more the lines \ref{P3-2D-appl-line03}--\ref{P3-2D-appl-line05}
in Algorithms~\ref{alg:MPPI:appl:2D} and~\ref{alg:MPPI:appl:3D}. Columns of $V_\ell$s, i.e.,
eigenvectors of $L_\ell$s are sines and cosines (see Lemmas~\ref{lem:spectrumP}--\ref{lem:spectrumM}).
Their outer products, i.e., the eigenvectors of $(-\FDL)$, are products of two or three goniometric
functions in two or three variables, respectively. Line \ref{P3-2D-appl-line03} can be then seen
as a kind of expansion of the updated residual into a Fourier series. After division of individual
components by the respective eigenvalues in line \ref{P3-2D-appl-line04}, the preconditioned residual
is obtained as a linear combination of these goniometric functions in line \ref{P3-2D-appl-line05}.
Thus, the pseudoinverse preconditioner is in our case essentially a Fourier solver.

Further, the Moore--Penrose pseudoinverse (i.e., in the regular case the inverse) preconditioner
can be used in principal as a complete stand-alone direct solver. However, since our matrices
might be singular (or close to singular), we need to introduce the threshold (see, e.g.,
line \ref{P3-2D-init-line05} in \ref{alg:MPPI:init:2D}) which complicates direct application
of such direct solver. Therefore, we decided to enclose it into outer loop of PCG. Such two-stage
algorithm is not significantly slower since most parts of the pseudoinverse is computed once.
On the other hand, one can anticipate, that such preconditioning causes rapid convergence of
the outer PCG.
\end{remark}

\section{Numerical results}\label{sec:numexp}

We present here several selected experiments we performed on data corresponding to three different
testing problems. More realistic results can be found in already published paper in Physical
Review B \cite{zigzagpaper} that is partially based on the presented approach.

\subsection{Testing problems}

The first problem tests basic scalability, it consists of four right-hand side matrices,
or better to say, four discretizations of the same continuous right-hand side. The second
problem contains only one right-hand side matrix and its purpose is to test different boundary
conditions. The third problem consists of four a bit more realistic right-hand sides --- it deals
with random distribution of charge and three of the right-hand sides are 3D. Thus, we present
here in total nine testing right-hand sides (six 2D and three 3D). All of these right-hand sides
are artificially fabricated by {\sc Ferrodo2}; see \cite{ferrodo}, \cite{zigzagpaper}.

\begin{figure}[b!]
\includegraphics[width=.49\textwidth]{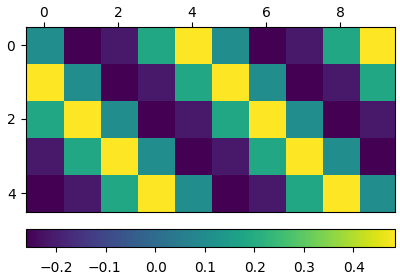}
\includegraphics[width=.49\textwidth]{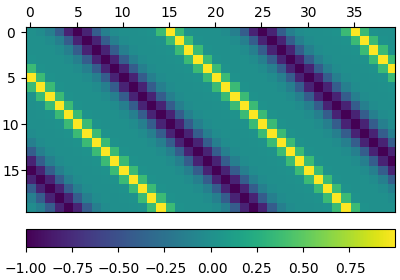}
\includegraphics[width=.49\textwidth]{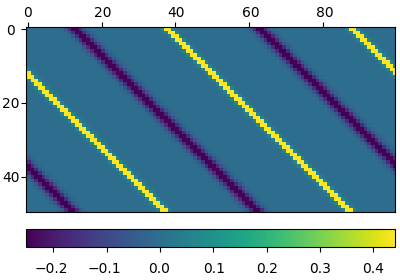}
\includegraphics[width=.49\textwidth]{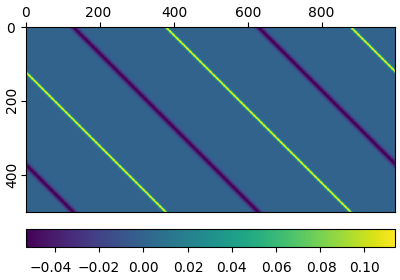}
\caption{Four right-hand sides of the regular two-periodic problem (see Problem~\ref{pb:1}).}\label{fig:pb:1}
\end{figure}

\begin{problem}[Regular two-periodic problem]\label{pb:1}
The first sample is a 2D problem, with regular alternating diagonal lines (or narrow stripes)
of positive and negative charges; see Figure~\ref{fig:pb:1}
(indices in the images are shifted by one compared to matrices, and start at $0$).
This problem is periodic in both physical dimensions, i.e., the corresponding Laplacian
is singular. We work with four discretizations of the same right-hand side, from the
coarsest to the finest: $5\times10$, $20\times40$, $50\times100$, and $500\times1000$.
\end{problem}

\begin{figure}[htb!]
\includegraphics[width=.49\textwidth]{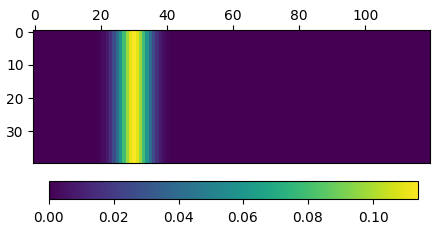}
\caption{Right-hand side of the problem with mixed boundary conditions (see Problem~\ref{pb:2}).}\label{fig:pb:2}
\end{figure}

\begin{problem}[Problem with mixed boundary conditions]\label{pb:2}
The second problem is also 2D, positive charge is localized in a single band; see Figure~\ref{fig:pb:2}.
On the left we prescribe zero electric potential, i.e., Dirichlet condition $u_\mathsf{L}=0$, whereas
on the right its derivative, i.e., the Newton condition $e_\mathsf{R}=-\frac12$. In the vertical
direction, the problem is periodic. We use only one discretization of moderate size $40\times 120$.
\end{problem}

\begin{figure}[htb!]
\includegraphics[width=.24\textwidth]{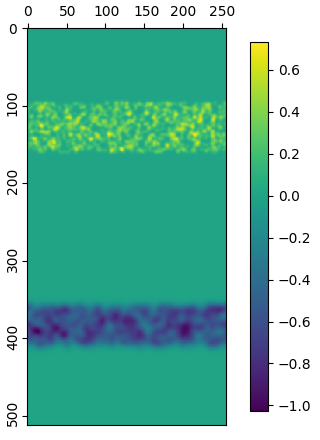}
\includegraphics[width=.24\textwidth]{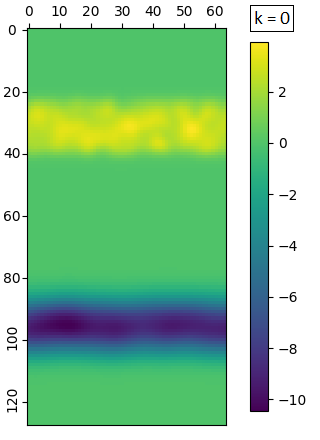}
\includegraphics[width=.24\textwidth]{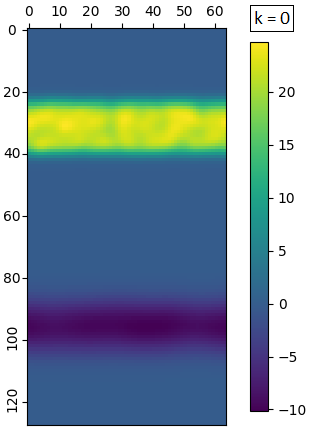}
\includegraphics[width=.24\textwidth]{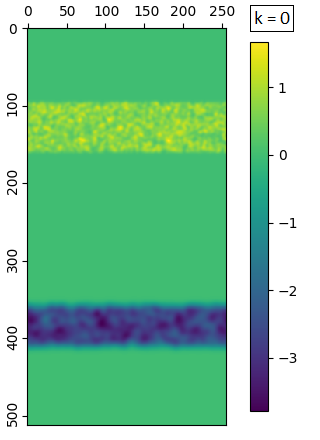}
\caption{Four right-hand sides of problem with randomly distributed charge (see Problem~\ref{pb:3}).
Notice that strips are of different width --- negative charge is more localized.
The last three images show the first frontal slices ($k=0$) of three-way tensors.}\label{fig:pb:3}
\end{figure}

\begin{problem}[Problem with randomly distributed localized charge]\label{pb:3}
The last problem is little bit more realistic. It is again periodic in all physical dimensions
(so the Laplacian is singular; similarly as in Problem \ref{pb:1}), the charge is localized in two
stripes, but randomly distributed; these stripes are not the same, the positive
charge spreads in wider strip than the negative.
This problem consists of four right-hand sides: the first
one is 2D $512\times 256$, the remaining three are 3D  $128\times64\times8$,
$128\times64\times64$, and $512\times256\times8$. Figure~\ref{fig:pb:3} shows the first
frontal slice in the case of 3D data (all slices look similarly).
\end{problem}

\subsection{Numerical experiments}

In all experiments we started the PCG with zero initial guess, i.e., $U_0=0_{n,q}$ 
in 2D, or $\mathcal{U}_0=0_{n,q,t}$ in 3D cases.  All the problems, in particular 
the right-hand sides (and thus also the initial residuals) are for convenience and 
better comparability of the results normalized such that 
\begin{equation}\label{eq:normalization}
  \|H\|_\fro = \|\vec(H)\|_2 = \frac1{nq} = \|R_0\|_\fro
  \quad \text{and} \quad
  \|\mathcal{H}\| \equiv \|\vec(\mathcal{H})\|_2 = \frac1{nqt} = \|\mathcal{R}_0\|,
\end{equation}
in the respective cases. In all the following plots we measure the effectiveness of 
individual methods (preconditioners) by the real run-time of the PCG. To obtain 
reasonable results:
\begin{itemize}\itemsep 0pt
\item First we experimentally fix the number of iterations for each pair of problem and preconditioner,
to get true residual small enough (we are not interested in stopping criterion at this moment; see 
Remark \ref{rmk:stoppincrit}).
\item Then we run the method several times (usually ten-times) and calculate the average run-time.
We measure the time for the initialization and the time for the iterative part separately.
\item Since the cost of individual iteration is independent on the iteration number, we
compute the time for one iteration simply by further averaging (the average time for the 
iterative part over total number of iterations).
\end{itemize}
These times are affected by collecting all the data durning the PCG run, however,
we collect the same data of the same size in each iteration for each preconditioner.
Presented results are computed by:
{\sc Matlab} R2024a (version 24.1.0.2653294, update~5)
running at
Microsoft Windows 10 Pro (version 10.0, build 19045, 64 bit),
Dell Latitude 5410 with Intel Core i5-10310U CPU @ 1.70/2.21 GHz and 8 GB RAM.

\begin{figure}[htb!]
\includegraphics[width=.8\textwidth]{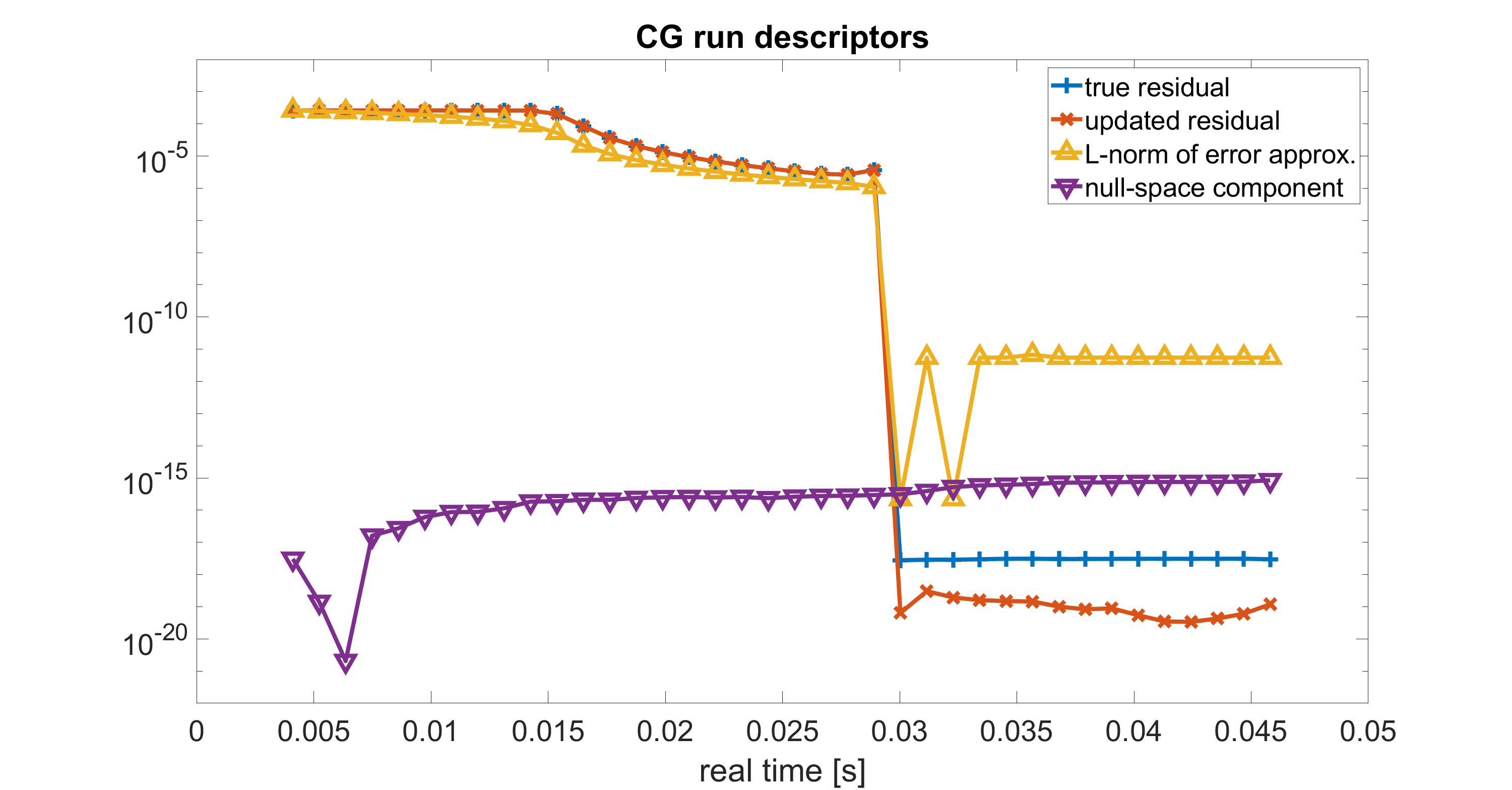}
\caption{Behavior of unpreconditioned CG observed by several quantities
(see Experiment~\ref{exp:1}).}\label{fig:exp:1-2}
\end{figure}

\begin{experiment}[Justification of method and convergence metric choices]\label{exp:1}
The first experiment is rather formalistic. We employ the second finest right-hand side of 
Problem \ref{pb:1}, i.e., we solve linear system with $5000$ unknowns. 

Since we already mentioned the Jacobi-like stationary iteration scheme in Section \ref{sssec:SSJLM},
we use it there for preconditioning, but we can use it as a stand-alone method. 
Figure \ref{fig:exp:1-1} compares convergence of the unpreconditioned CG and Jacobi-like
method for three different values of $\omega$. One can easily see the dominance of CG.
Notice that the (horizontal) time axis is in logarithmic scale, due to disproportional
run-times of both methods. Notice also that we need more than $6500$ of Jacobi iterations
to converge, thus we mark only each $250$th.

\begin{figure}[htb!]
\includegraphics[width=.8\textwidth]{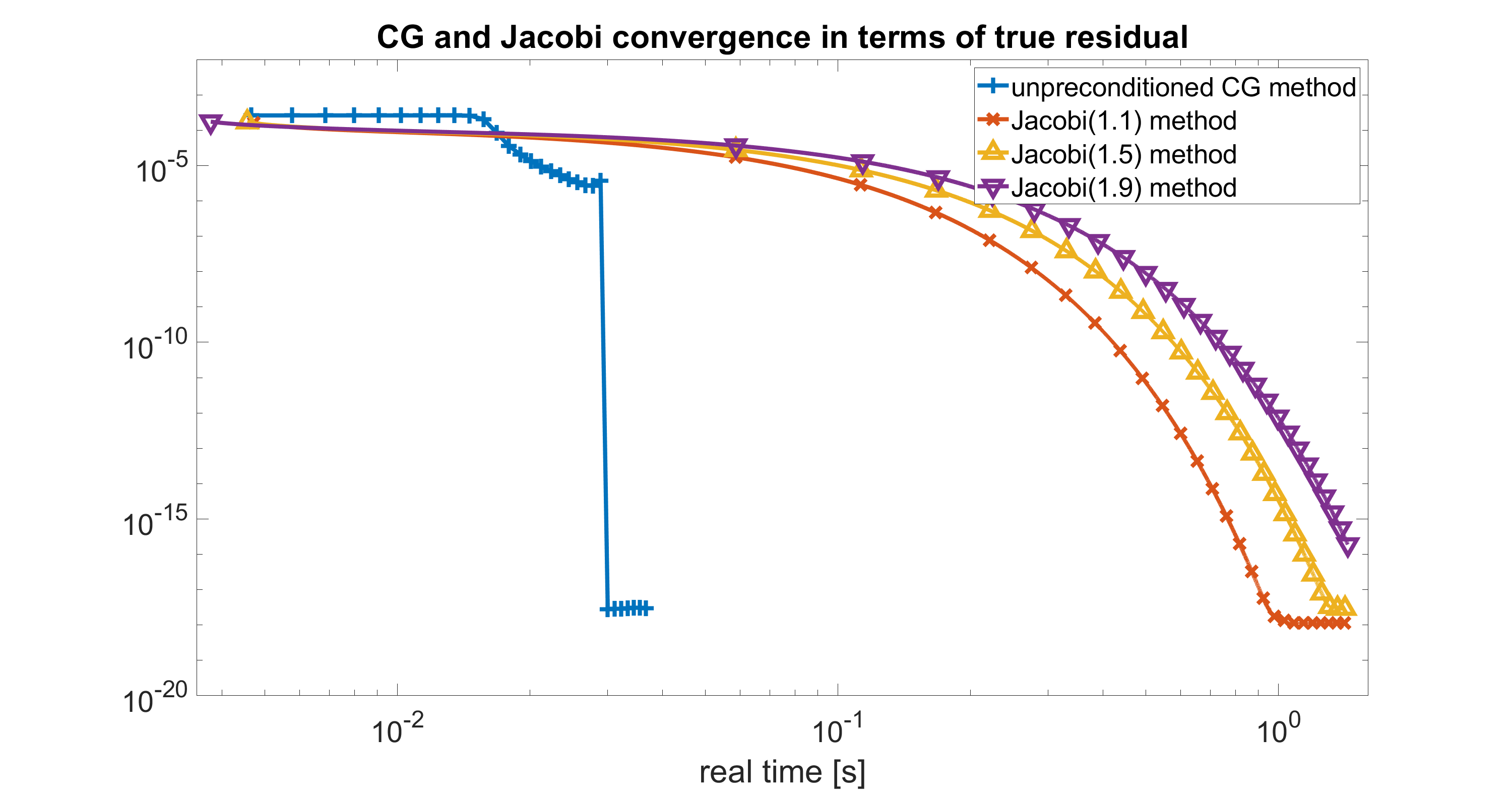}
\caption{Convergence of unpreconditioned CG and Jacobi-like
stationary iterative method for different values of $\omega$. 
Time axis is in logarithmic scale; markers at convergence 
curves for Jacobi denotes each $250$th iterations
(see Experiment~\ref{exp:1}).}\label{fig:exp:1-1}
\end{figure}

We measure the convergence by the norm of the true residual (i.e., $\|h-Lu_s\|_2$ in 1D case).
This choice is definitely wrong for CG/PCG (see also Remark \ref{rmk:stoppincrit}). However in 
our relatively simple problem of Poisson's equation, it seems to be applicable. 
Figure \ref{fig:exp:1-2} shows behavior of the true residual and the computed residual 
(i.e., norm $\sqrt{\rho_s}=\|r_s\|$ from Algorithm \ref{alg:PCG:3D}) --- both quantities are
indistinguishable up to the moment of convergence, then the true residual stagnates while
computed further decrease; see, e.g., \cite{b:meurant06} or \cite{b:liesenstrakos12}. We also ploted quantity
that behaves like the $L$-norm of error. By the $L$-norm of error we understand (again in 1D case, for simplicity) 
\[
  \|u_s-u\|_L 
  = \sqrt{(u_s-u)^\trp L(u_s-u)}
  = \sqrt{u_s^\trp Lu_s - 2u_s^\trp Lu + u^\trp Lu}\geq0.
\] 
Since $u$ is the exact solution, then $Lu = h$. Thus
\[ 
  \|u_s-u\|_L^2  = \underbrace{u_s^\trp Lu_s - 2u_s^\trp h}_{\displaystyle\kappa_s} \,+\, \|u\|_L^2  \geq 0,
\]
where the underbraced part $\kappa_s$ is fully computable, and (the square of) the $L$-norm 
of $u$ is unknown, but constant. Thus we compute all $\kappa_s$s, 
shift them to get nonnegative quantities, e.g., $\kappa_s-\min_s(\kappa_s)$, 
and square-root them; altogether we approximate 
\[
  \|u_s-u\|_L \approx \eta_s \equiv \sqrt{\kappa_s-\min_s(\kappa_s)} \geq 0. 
\]
The final quantity $\eta_s$ is in the plot further scaled and shifted for better
comparability with the residuals; we plot in particular
\[
  \alpha\eta_s + \beta, 
  \qquad\text{where}\qquad
  \alpha = \frac{\|h-Lu_1\|_2}{\eta_1}, 
  \quad\text{and}\quad
  \beta = \varepsilon_{\mathsf{M}} = 2^{-52} \approx 2.2204\,\cdot10^{-16}
\]
is the machine precision. Our goal not to get robust approximation, but to observe 
possible drops in the norm. Figure \ref{fig:exp:1-2} shows that the most significant 
drop in all of the three mentioned descriptors (both residuals and $L$-norm of error 
approximation) coincides. Notice that we complement this information 
by  fourth descriptor: Since we use problem with singular matrix, we also observe
the norm of the null-space component in the solution. In all the remaining plots
we restrict ourselves only to the true residual.
\end{experiment}

\begin{experiment}[Performance of various preconditioners]\label{exp:2}
In the second experiment we again use the second finest right-hand side of Problem \ref{pb:1}. 
We compare PCG with individual preconditioners, each with some basic settings. In particular:
$\mathfrak{M}=I$ (no preconditioner),
$\mathfrak{M}_{\text{\rm Jacobi}(3,1.3)}$ (three steps of Jacobi-like method, i.e., $p=3$, with $\omega=1.3$),
$\mathfrak{M}_{\text{\rm low-rank}(3)}$ (Kronecker rank $r_\kro = 3$ approximate inverse), and
$\mathfrak{M}_{\text{\rm p-inv}}$ (pseudoinverse preconditioner).

\begin{figure}[htb!]
\includegraphics[width=.8\textwidth]{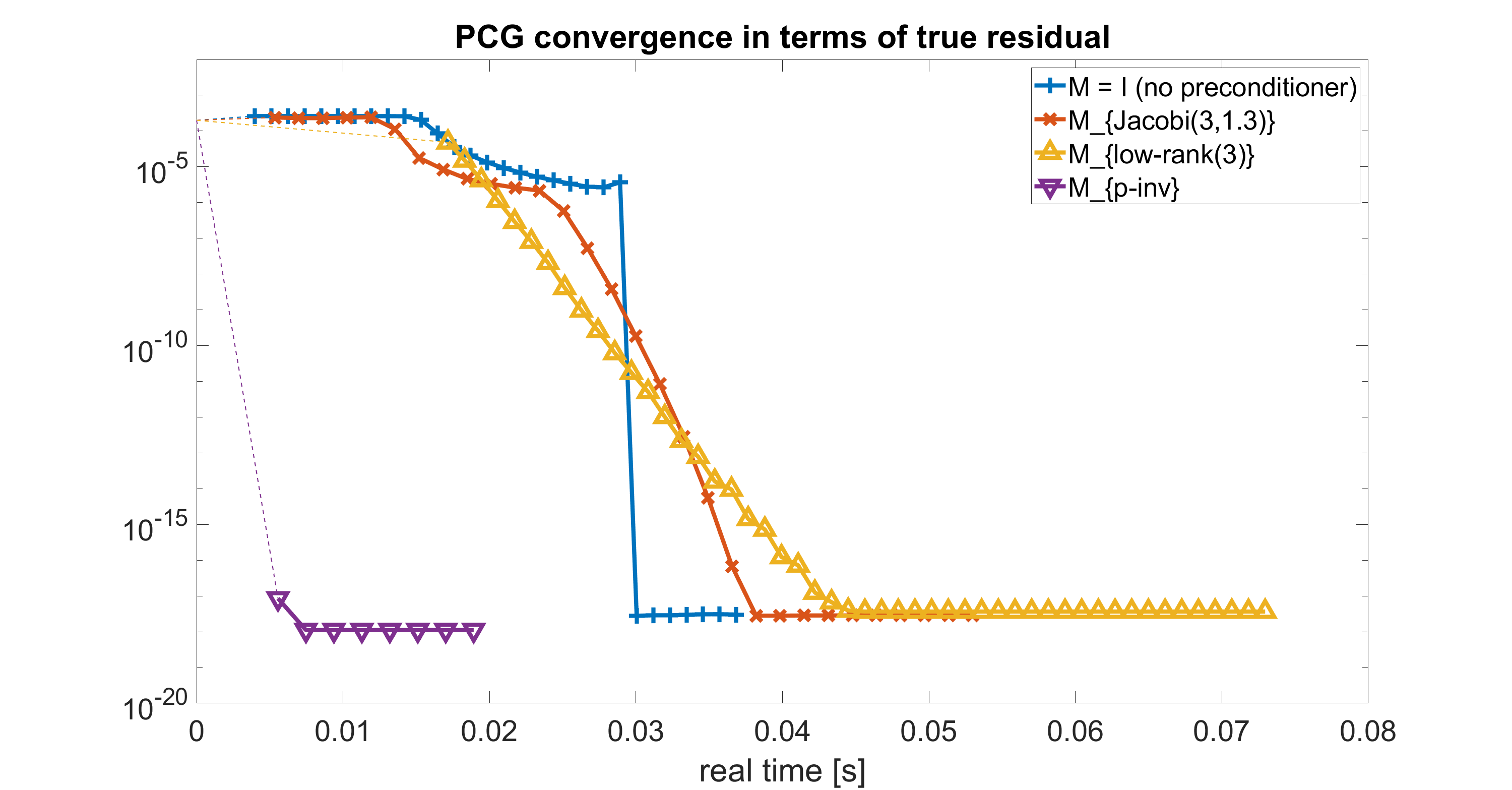}
\caption{PCG with various preconditioners applied on Problem \ref{pb:1}
of size $50\times100$ (see Experiment~\ref{exp:2}).}\label{fig:exp:2-1}
\end{figure}

We see that the unpreconditioned CG starts fastest as expected. On the other hand,
the slowest start has the PCG with low-Kronecker-rank preconditioner, as already
discussed in the last paragraph of Section \ref{sssec:AILKR}. Its convergence is
very straightforward, however slower than for unpreconditioned CG. Nevertheless,
PCG with the Moore--Penrose pseudoinverse-based preconditioner beats them all ---
its convergence is very rapid. After the first iteration it is already almost done
(we use the almost inverse for preconditioning), however the second iteration further
pushes down the true residual.

For completeness, Figures \ref{fig:exp:2-2} and \ref{fig:exp:2-3} illustrate behavior
of parametrized preconditioners $\mathfrak{M}_{\text{\rm Jacobi}(p,\omega)}$
and $\mathfrak{M}_{\text{\rm low-rank}(r_\kro)}$, respectively. One can see that
the low-Kronecker-rank preconditioner is getting closer to the pseudoinverse 
preconditioner with higher rank, but with more expensive initialization. The artifact
on the already converged curve for $r_\kro$ is probably due to the indefiniteness
of the preconditioner (which was indicated for $r_\kro = 4$ and $7$).  
\end{experiment}

\begin{figure}[htb!]
\includegraphics[width=.8\textwidth]{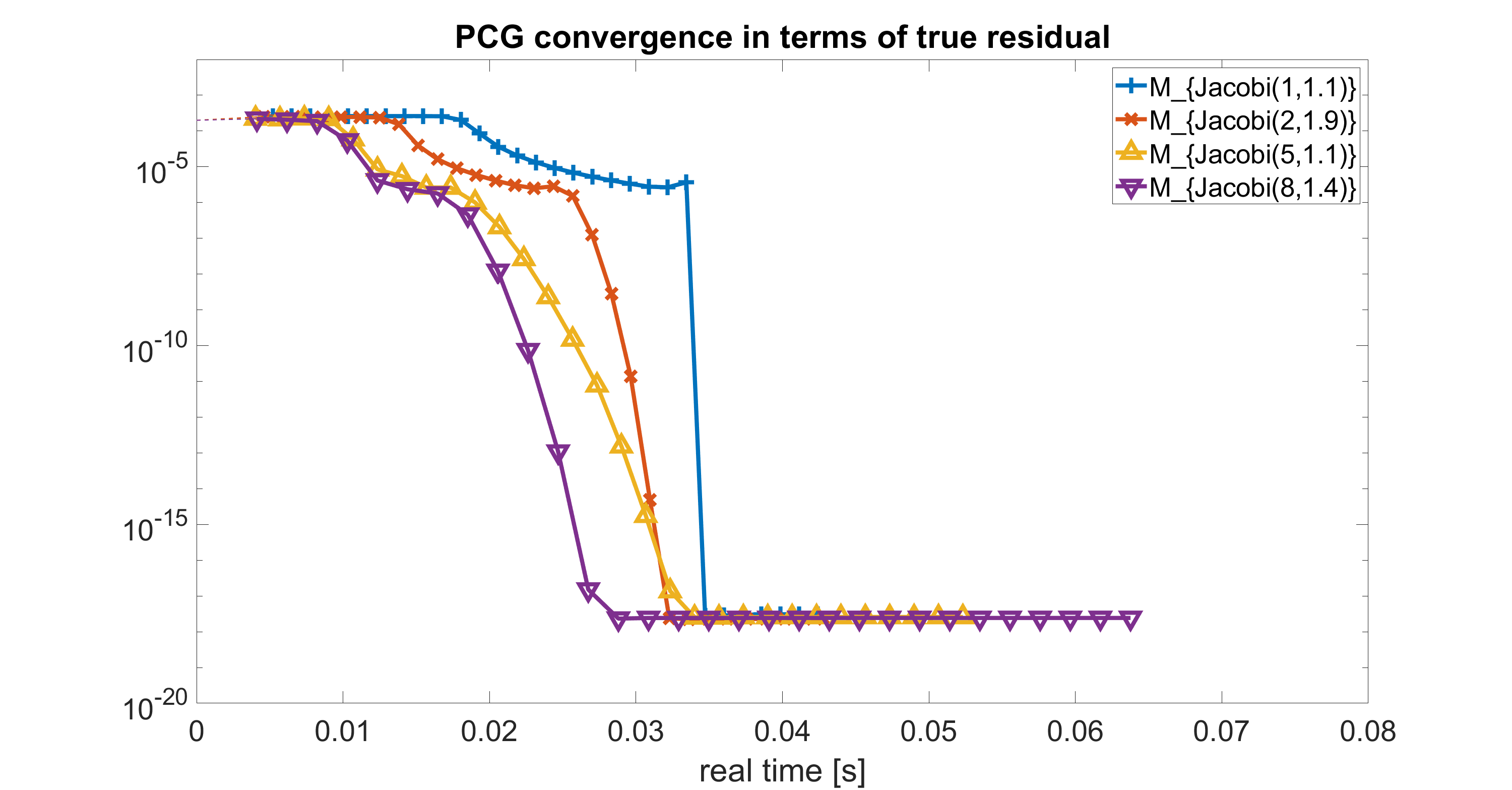}
\caption{Performance of Jacobi-like preconditioner is only weakly involved
by the number of inner iterations $p$ or the value of the splitting parameter $\omega$
(compare with Figure~\ref{fig:exp:2-1}).}\label{fig:exp:2-2}
\end{figure}

\begin{figure}[htb!]
\includegraphics[width=.8\textwidth]{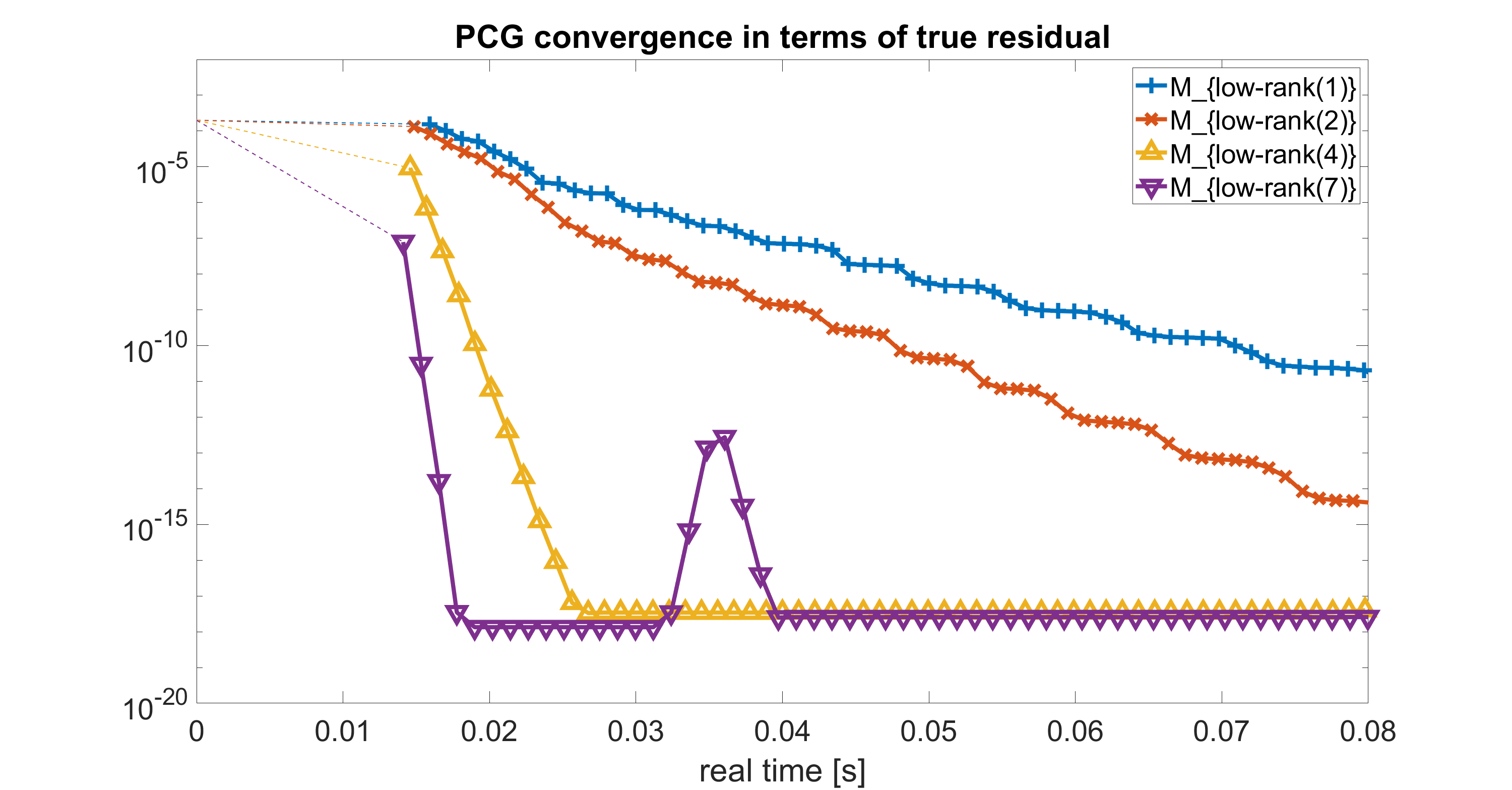}
\caption{Low-Kronecker-rank approximation of inverse works as preconditioner
significantly better with higher ranks $r_\kro$ (if it works, see the last paragraph
of Section \ref{sssec:AILKR}; compare with Figure~\ref{fig:exp:2-1}).}\label{fig:exp:2-3}
\end{figure}

\begin{experiment}\label{exp:3}
The third experiment is to verify the performance of the chosen preconditioner --- the
pseudoinverse-based --- on all three testing problems, all nine right-hand sides.
For the results see Figure \ref{fig:exp:3}.

\begin{figure}[htb!]
\includegraphics[width=.8\textwidth]{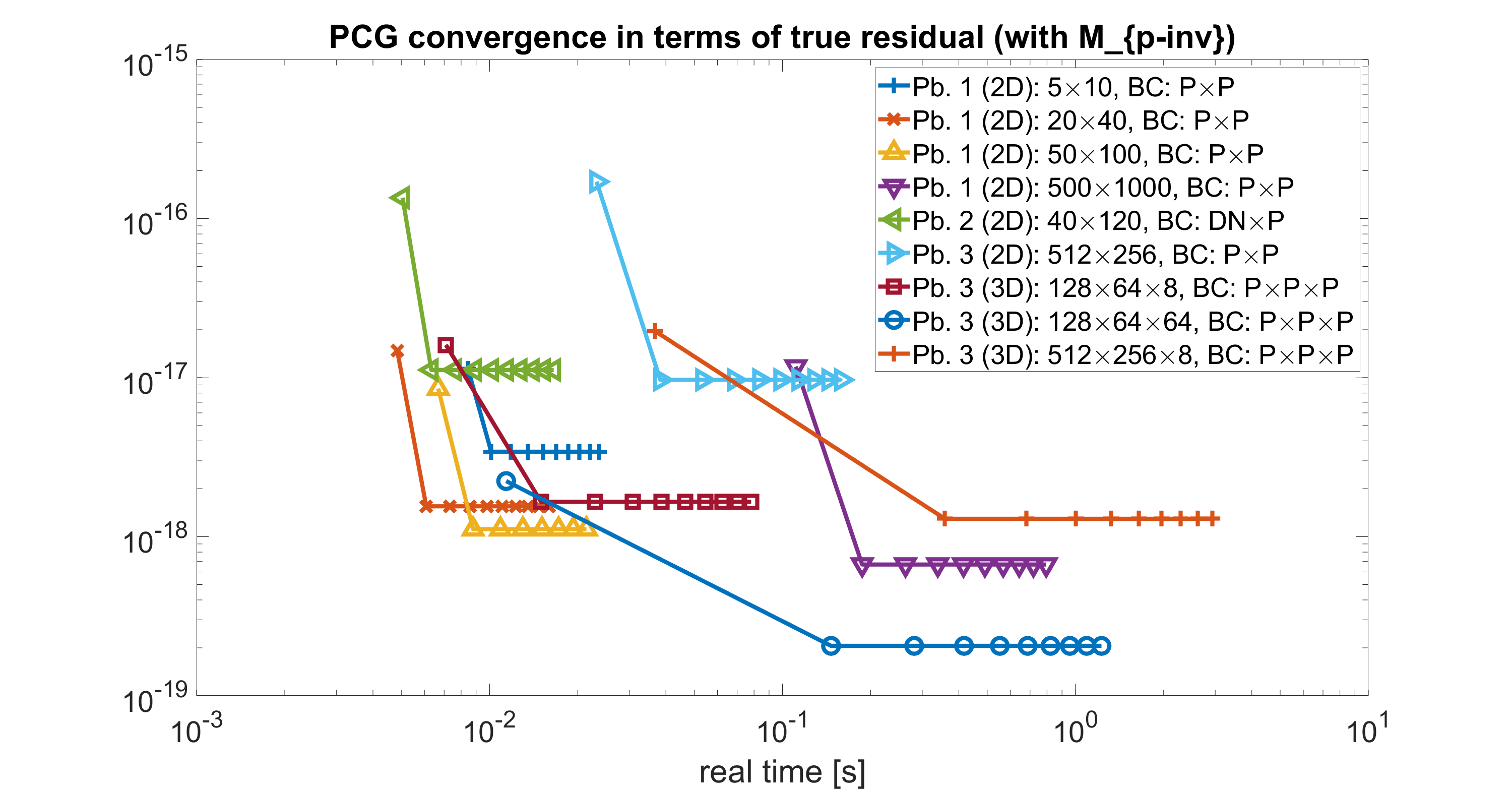}
\caption{The convergence curve of PCG with the pseudoinverse-based preconditioner
is almost problem independent. Time axis is in logarithmic scale 
(see Experiment~\ref{exp:3}).}\label{fig:exp:3}
\end{figure}

Notice that the first and smallest system is of $50$ equations, whereas the 
largest is of $1\,048\,576$ equations. PCG for such systems converge on
disproportional time-scales. To make the results visible in one plot, the 
time axis needs to be also in logarithmic scale. Therefore, we cannot
start the curve at zero time; thus the convergence curves start at the time
after initialization and show only the iterative part. This brings another problem:
Since the convergence is very rapid, already in the first iteration, and initial part
is missing, we cannot directly see from where it converges whereto. Recall that 
the norm of the initial residual is normalized (\ref{eq:normalization}); i.e.,
$\|R_0\|_\fro = 2 \,\cdot 10^{-2}$ for the smallest, and 
$\|\mathcal{R}_0\| = 9.5367431640625 \,\cdot 10^{-7}$ for the largest system.
\end{experiment}

\subsection{Solutions}
\begin{figure}[htb!]
\includegraphics[width=.8\textwidth]{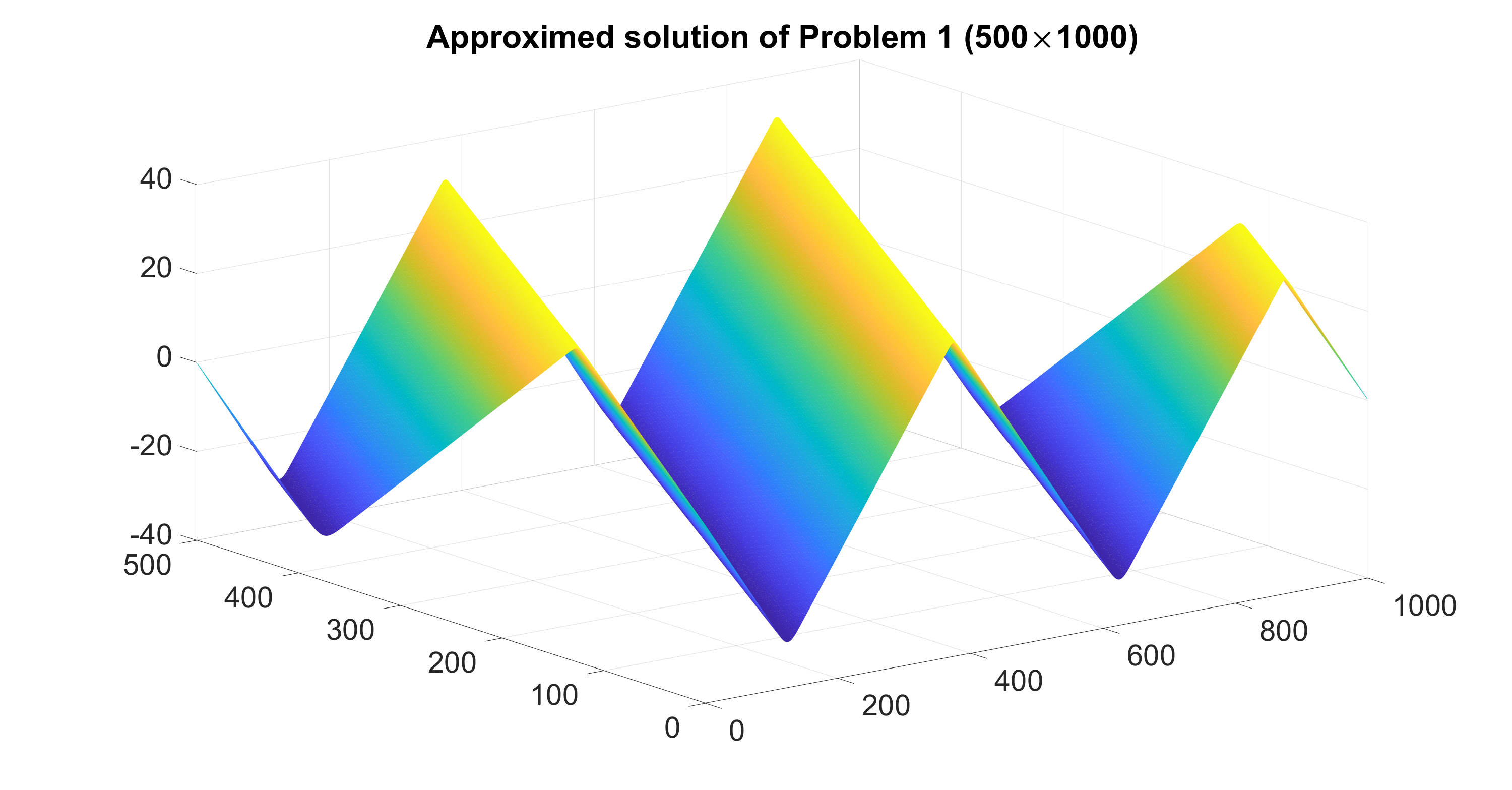}
\includegraphics[width=.8\textwidth]{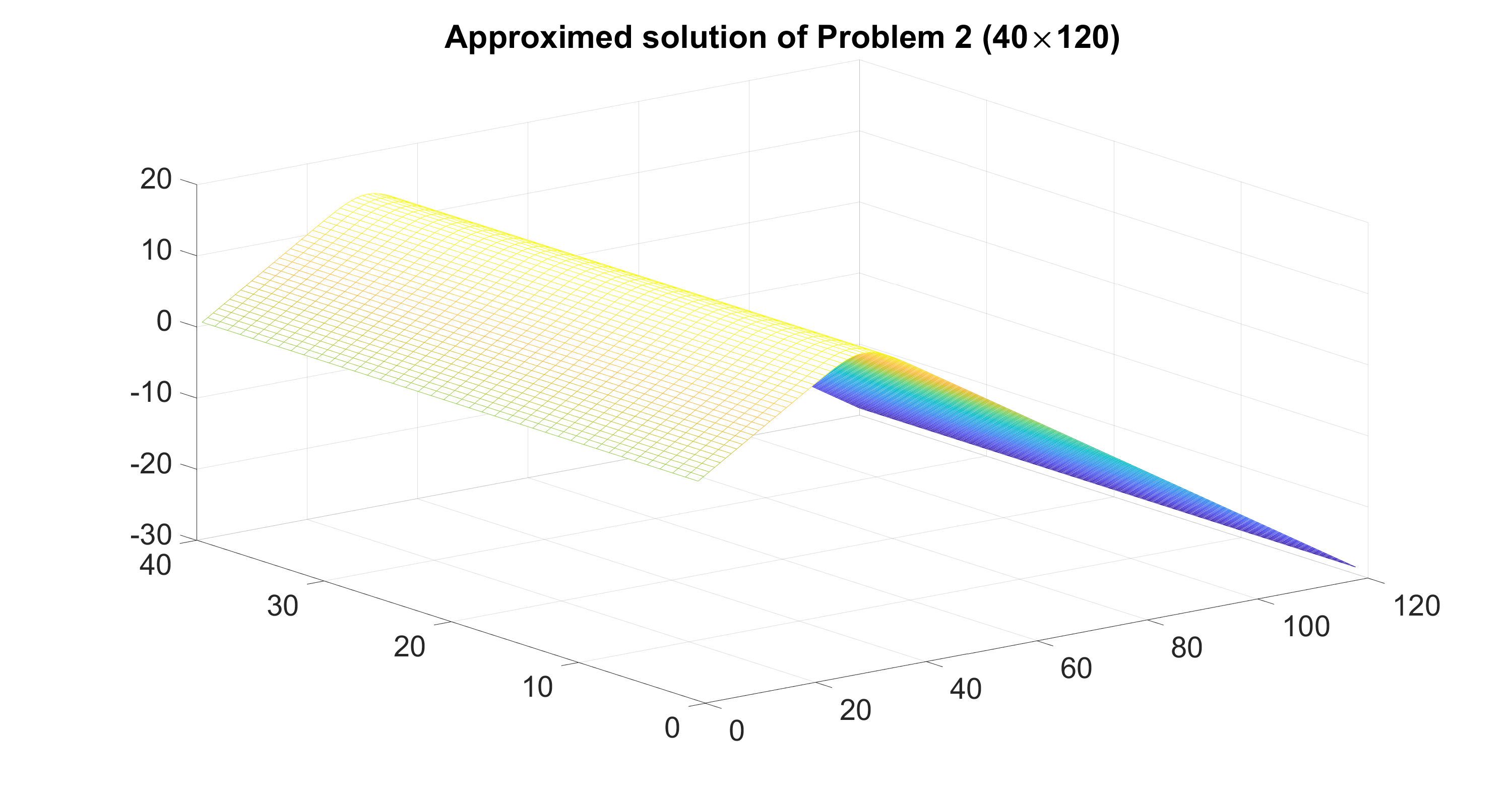}
\includegraphics[width=.8\textwidth]{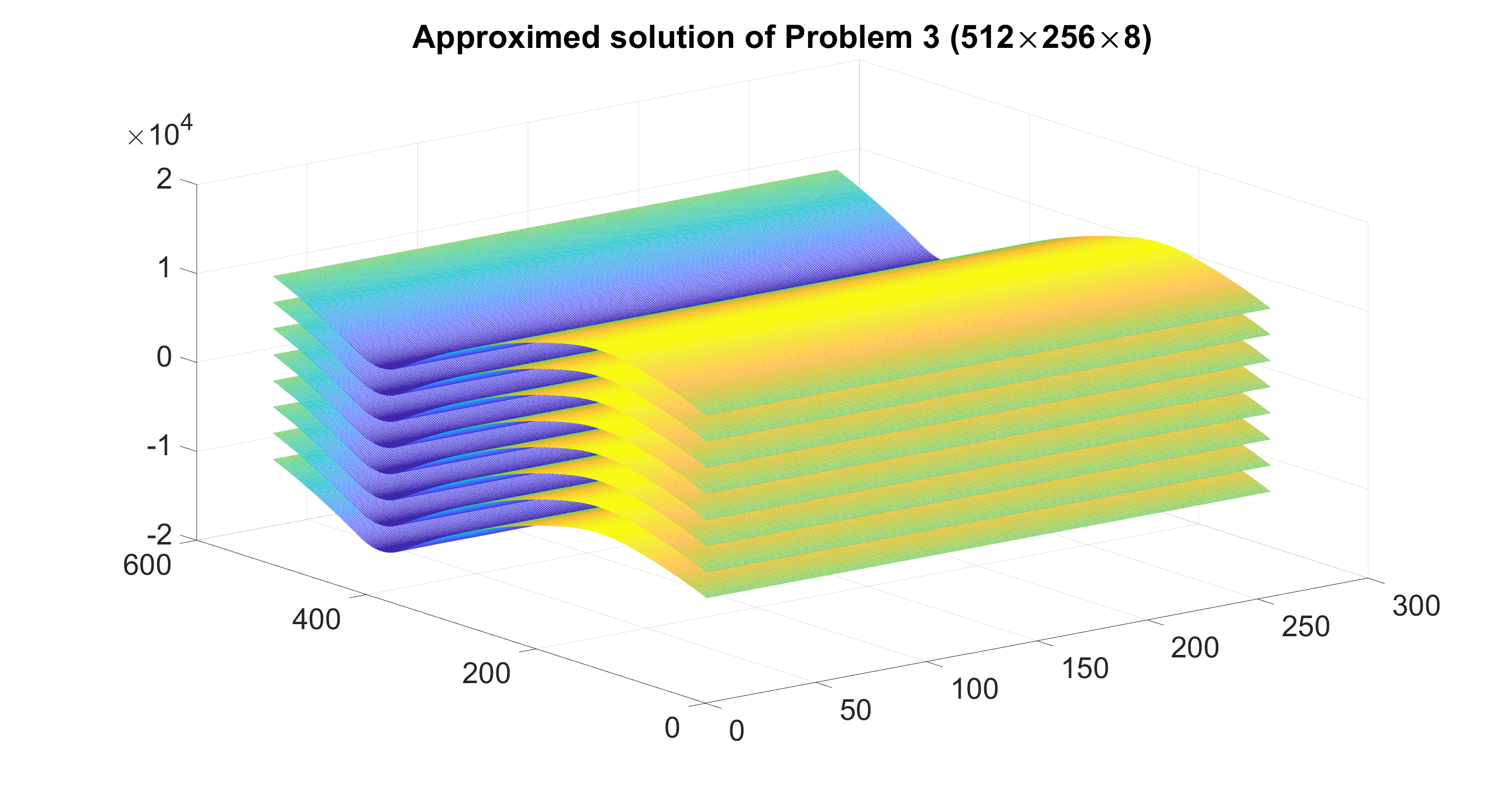}
\caption{Approximate solutions of all three problems. Vertical axis 
measures the electric potential (in the bottom plot, there is the
vertical scale related to the fifth layer from the top).}\label{fig:solutions}
\end{figure}

Just for completeness we present the approximate solutions of all three problems, 
computed by PCG with pseudoinverse-based preconditioner, in Figure~\ref{fig:solutions}. 
Since in Problem \ref{pb:1} (the top plot) the charge is localized in thin parallel 
and equidistant lines, the electric potential changes linearly between these lines.
In the solution of Problem \ref{pb:2} (the middle plot) the electric potential
is fixed on the left edge at zero and on the right it falls down with the rate $-1/2$,
as prescribed by the boundary conditions. 
The unequal width of two stripes with randomly distributed charges in Problem \ref{pb:3}
(the bottom plot) causes different curvature of the electric potential.

\section{Conclusions}\label{sec:conclusions}

We have presented an implementation of the preconditioned method of conjugate gradient (CG)
that is suitable for Kronecker-structured problems. In our case particularly the
Poisson's problem originated in ferroelectricity. We have briefly analyzed several 
preconditioners, their behavior and also computational cost in comparison to the CG itself.

Since the right-hand side is dense with general rank(s) (thus we cannot exploit 
its low-rank structure and employ the low-rank arithmetic) and of shape 
close enough to square (or cube, or hyper-cube), we are in principle able to 
store also the individual objects of the Kronecker-structured  system operator
in dense format. Consequently, we were able to reasonably well approximate the 
Moore--Penrose pseudoinverse of the operator, and use this approximation for 
preconditioning of CG. 

This way preconditioned CG converges rapidly. So rapidly that only a few iterations
is necessary --- in our case of a simple Laplace operator, it seems that two or three 
are  sufficient (but we perform ten, to be sure; see Figure~\ref{fig:exp:3}). 
Perhaps, instead of a {\em CG with pseudoinverse-based preconditioner}, it would be better 
to call this approach a {\em pseudoinverse-based direct solver with refinement by CG-iteration}.

The presented implementation is done in {\sc Matlab} and serves mostly for testing of methods 
and  preconditioners. We prepared also C/C++ version that employs the Intel oneAPI Math Kernel 
Library (MKL); see \cite{MKL1}, \cite{MKL2}. This version is implemented into, and 
practically used  in the code {\sc Ferrodo2}; see \cite{ferrodo}, \cite{zigzagpaper}.
First results computed by this code have been already published in Physical
Review B \cite{zigzagpaper}. 

\section*{Acknowledgements}

We wish to thank to Jana \v{Z}\'{a}kov\'{a} for her help with preparation of the manuscript
of this paper.

{\small
{\em Authors' addresses}:
{\em V\v{e}nceslav Chumchal},
  Department of Mathematics, Technical University of Liberec, Studentsk\'{a} 1402/2, 461 17 Liberec 1, Czech Republic,
  e-mail: \texttt{venceslav.chumchal@\allowbreak tul.cz} \& \texttt{venceslav.chumchal@\allowbreak gmail.com}.
{\em Pavel Marton},
  Institute of Physics, Czech Academy of Sciences, Na Slovance 2, 182 00 Praha 8, Czech Republic,
  \& Institute of Mechatronics and Computer Engineering, Technical University of Liberec, Studentsk\'{a} 1402/2, 461 17 Liberec, Czech Republic,
  e-mail: \texttt{pavel.marton@\allowbreak tul.cz}.
{\em Martin Ple\v{s}inger} (corresponding author),
  Department of Mathematics, Technical University of Liberec, Studentsk\'{a} 1402/2, 461 17 Liberec 1, Czech Republic,
  e-mail: \texttt{martin.plesinger@\allowbreak tul.cz}.
{\em Martina \v{S}im\r{u}nkov\'{a}},
  Department of Mathematics, Technical University of Liberec, Studentsk\'{a} 1402/2, 461 17 Liberec 1, Czech Republic,
  e-mail: \texttt{martina.simunkova@\allowbreak tul.cz}.
}

\end{document}